\documentstyle[aps,amstex,pre,multirow,preprint,lscape,epsfig,subfigure,hhline]{revtex}

\def\f{\frac}
\def\l{\left}
\def\r{\right}
\def\p{\partial}

\begin{document}


\title{High resolution conjugate filters 
        for the simulation of flows}

\author{Y. C. Zhou  and  G. W. Wei\footnote{Corresponding author}}

\address{Department of Computational Science, 
     National University of Singapore,
     Singapore 117543} 

\date{\today}
\maketitle

\begin{abstract}

This paper proposes a Hermite-kernel realization of the 
conjugate filter oscillation reduction (CFOR) scheme for 
the simulation of fluid flows. The Hermite kernel is 
constructed by using the discrete singular convolution (DSC)
algorithm, which provides a systematic generation of 
low-pass filter and its conjugate high-pass filters.
The high-pass filters are utilized for approximating spatial 
derivatives in solving flow equations, while the conjugate 
low-pass filter is activated to eliminate spurious oscillations 
accumulated during the time evolution of a flow. As both low-pass 
and high-pass filters are derived from the Hermite kernel, 
they have similar regularity, time-frequency localization, 
effective frequency band and compact support. Fourier analysis 
indicates that the CFOR-Hermite scheme yields a nearly optimal 
resolution and has a better approximation to the ideal low-pass 
filter than previously CFOR schemes. Thus, it has better potential 
for resolving natural high frequency oscillations from a shock.
Extensive one- and two-dimensional numerical examples, including 
both incompressible and compressible flows, with or without 
shocks, are employed to explore the utility, test the resolution, 
and examine the stability of the present CFOR-Hermite scheme. 
Extremely small ratio of point-per-wavelength (PPW) is achieved 
in solving the Taylor problem, advancing a  wavepacket 
and resolving a shock/entropy wave interaction. The present 
results for the advection of an isentropic vortex compare very 
favorably to those in the literature.


\end{abstract}

\section{Introduction}

The rapid progress in constructing high performance computers 
gives great impetus to the development of innovative high order 
computational methods for solving problems in computational 
fluid dynamics (CFD). With high-order numerical methods, complex 
flow structures are expected to be resolved more accurately in both  
space and time. For flow problems with sophisticate structures,  
high resolution is a must in order for the structural information
to be correctly extracted. High order numerical resolution is of 
crucial importance in all kinds of simulations of turbulence. For 
example, direct numerical simulation (DNS) of turbulence is one of 
the typical situations  where high resolution schemes are definitely 
required because the amplitude of the Fourier response of the 
velocity field is continuously distributed over a wide range of 
wavenumbers.  The situation is more complex in the large eddy 
simulation (LES), where the energy spectrum is not only continuous,
but also slow decaying at sufficiently high wavenumbers. Unlike 
a linear problem, the numerical error of the $q$th order finite 
difference scheme is not linearly scaled with the power of grid
size $\Delta^q$. Thus numerical errors have to be strictly controlled 
by choosing a grid which is much smaller than the smallest spatial 
structure in turbulent motion\cite{ghosal96}. The simulation become 
more difficult as the Reynolds number increases. For incompressible 
Euler flows, although the total kinetic energy is conserved, 
its spectral distribution  of the solution  could move towards 
the high frequency end with respect to the time evolution.
As a consequence, the integration will collapse at a sufficiently 
long  time with a finite computational grid due to the limited high 
frequency resolution  of a given scheme. Nevertheless, a higher 
resolution scheme can dramatically improve the stability of the 
numerical iterations compared with a lower resolution one because 
of its ability to better resolve the high frequency components. Using 
three candidate schemes, i.e., a second-order finite difference 
scheme, a forth-order finite difference scheme and a pseudo`spectral 
scheme,  Browning and Kreiss\cite{browning89} demonstrated that 
in the long-time  DNS of turbulence, only the forth-order 
scheme can yield results comparable with those of the pseudospectral 
scheme.

Apart from the DNS and LES of turbulence, computational aeroacoustics 
(CAA) is another area where high resolution schemes are highly 
demanded  for the simulation of waves with large spectral bandwidth 
and amplitude disparity\cite{tam95}. For an imperfectly expanded 
supersonic jet, the Strouhal number ranges from $10^{-2}$ to $10^{1}$. 
The velocity fluctuation of the radiated sound wave can be four 
orders of magnitude smaller than that of the mean flow. Therefore,
the real space and time numerical simulation of turbulence 
flow and sound wave interaction in an open  area and  long 
time process is still a severe challenge to practitioners of CAA.
Obviously, high resolution schemes can be used to alleviate
the demanding on a dense mesh in this situation.

The concept of high resolution is subjective. For a given 
scheme, it might provide a relatively high resolution for one 
problem but fail doing so for another problem. Yang et 
al\cite{yang01} show that for the heat equation, a local scheme 
delivers higher accuracy than the Fourier pseudospectral method, 
while the latter achieves  better accuracy in solving 
the wave equation. Essentially, the pseudospectral method
possesses the highest numerical resolution for approximating 
bandlimited periodic $L^2$ functions. However, it might not 
be the most accurate method for other problems. It is noted 
that much argument given to the numerical resolution of 
computational schemes in the literature is analyzed with 
respect to the discrete Fourier transform. Such a resolution 
should be called {\em Fourier resolution} and differs much 
from the numerical resolution in general. It is commonly 
believed that a scheme of high Fourier resolution, i.e.,
it  provides a good approximation to the Fourier transform 
of the derivative operators over a wider range of wavenumbers, 
will perform well in numerical computations. Although the 
results of discrete Fourier analysis might be consistent with 
the numerical resolution and be useful for a large  class 
of problems, they are strictly valid only for bandlimited 
periodic $L^2$ functions. For example, the Fourier resolution 
of the standard finite difference schemes is not very high 
as given by the Fourier analysis. However, the standard finite 
difference schemes are the exact schemes for approximating 
appropriate polynomials.

There are two approaches to obtain high resolutions in a 
numerical computation. One is to employ a spectral method,
such as the pseudospectral or Chebyshev method. In general
spectral methods provide very high numerical resolution for a 
wide variety of physical problems. However, they often have 
stringent constraints on applicable boundary conditions and 
geometries. The other approach is to modified the coefficients
of the standard finite difference scheme so that the high
frequency components of a function are better approximated 
under the cost of the approximation accuracy for the low 
wavenumbers. Typical examples include the compact scheme 
\cite{lele92} and dispersion relation preserving (DRP) 
scheme \cite{tam93}. Fourier analysis of the both schemes 
indicates that they can give a better representation in a 
broader wavenumber range than their central difference 
counterparts. Therefore, they might yield a good resolution 
for small flow structures in a more stable manner.  These 
modified finite difference schemes are very popular and give 
better numerical results for many physical problems.

The presence of shock waves adds an extra level of difficulty
to the seeking of high resolution solutions. Formally, a shock 
will lead to a first-order error which will propagate to the 
region away from the discontinuity in a solution obtained by using 
a high order method\cite{engquist98}. Therefore, the attainable
overall order of resolution is limited. However, for a given 
problem, the approximation accuracy achieved by using a high 
resolution scheme can be much higher than that obtained by using 
a low resolution scheme. In the context of shock-capturing, high 
resolution  refers this property not only in smooth regions but 
also in regions  as close to the discontinuity as possible.
One of popular high order shock-capturing schemes is an essential 
non-oscillatory (ENO) scheme\cite{HEOC,Shu1,EShu}. The key
idea of the ENO scheme is to suppress spurious oscillations  near 
the shock or discontinuities, while maintaining a  higher-order 
accuracy at smooth regions. This approach was  further 
extended into  a weighted  essentially  non-oscillatory 
(WENO) scheme\cite{LOC}.  The WENO approach takes a linear 
combination of a number of high-order  schemes of both central 
difference  and  up-wind type.  The central difference type of
schemes has a larger weight factor at the  smooth region while 
the up-wind type of schemes plays a major role near the shock  or 
discontinuity. The arbitrarily high order schemes which utilize 
the hyperbolic Riemann problem for the advection of the higher 
order derivatives (ADER)\cite{toro98} and monotonicity preserving 
WENO scheme\cite{balsara00} are the most recent attempts in
developing higher-order shock-capturing schemes. These schemes 
are very efficient for many hyperbolic conservation law systems.

One of the most difficult situations in the CFD is the interaction 
of shock waves and turbulence, which occurs commonly in high speed 
flows such as buffeting, air intaking and jet exhaust.  Many high 
order shock-capturing schemes are found to be too dissipative to be 
applicable for the long-time integration of such flows\cite{LLM,GMSCD}. 
In the framework of free decaying turbulence, the effect of a
subgrid-scale model was masked by some high-order shock-capturing 
schemes. Excessive dissipation degrades the numerical resolution 
and smears small-scale structures. Typically, there are two ways to reduce
the excessive numerical dissipation. One is to rely on shock 
sensors for switching on-off dissipative terms locally. The 
spatial localization of the sensors is crucial to the success 
of the approach. In the framework of synchronization, a set of 
interesting nonlinear sensors was proposed\cite{weiprl} for 
achieving optimal localization automatically. The other approach 
is to develop low dissipation schemes, such as filter schemes.  
In fact, a filter, or a dissipative term, can be extracted from 
a regular scheme, such as the ENO. Such a term can be  combined 
with another high-order non-dissipative scheme, such as a central 
difference scheme, to  construct a new shock-capturing scheme 
\cite{yee99,eric01}. The overall dissipation of the new scheme 
can be controlled by using an appropriate switch. Numerical 
experiments have shown that these schemes can be of higher 
accuracy and less dissipative than the original shock-capturing 
schemes\cite{eric01}.

Recently, a conjugate filter oscillation reduction (CFOR) has 
been proposed as a nearly optimal approach for treating hyperbolic 
conservation laws\cite{weigu,guwei}. The CFOR seeks a uniform high 
resolution both in approximating derivatives and in filtering 
spurious oscillations. The conjugated filters are constructed by 
using the discrete singular convolution (DSC) algorithm, which 
provides a spectral-like resolution in a smooth region. Moreover, 
the accuracy of the DSC algorithm is controllable by adjusting the 
width of the computational support. From the point of signal 
processing, the DSC approximations of derivatives are high-pass 
filters, whereas, spurious oscillations are removed by using 
a conjugate low-pass filter, without resort to upwind type of 
dissipative terms. The  low-pass and high-pass filters are 
conjugated in the sense that they are generated from the same 
DSC kernel expression. The resolution of the CFOR scheme is 
therefore determined by the common effective frequency band of 
the high-pass and low-pass filters. The effective frequency band 
is the wavenumber region on which both low-pass and high-pass  
filters have correct frequency responses. Although, the errors 
produced at a discontinuity are still of the first order, their 
absolute values are very small due the use of the high 
resolution CFOR scheme.

There are a variety of DSC kernels that can be used for realizing 
the CFOR scheme\cite{weijcp}. In fact, these kernels might exhibit 
very different characteristic in numerical computations, in 
particular, for shock-capturing. The latter is very sensitive to 
the property of DSC kernels. So far, the most widely used kernels 
is the regularized Shannon kernel (RSK)\cite{weigu,guwei}. The 
success of the RSK in numerical computations, particularly, in 
shock capturing, motivates us to search for new DSC kernels 
which might have more desirable features for certain aspects of 
shock-capturing. Indeed, Hermite kernel (HK) is a very special
candidate\cite{Korevaar,Hoffman,weijpa}. Compared with the RSK, 
the Hermite kernel has a narrow effective frequency band but its 
low-pass filter better approximates the idea low-pass filter. As 
such, we speculate that the Hermite kernel should have a resolution 
similar to the RSK in the smooth flow region and might have better 
overall performance in the presence of shock. Extensive numerical 
tests are given in this paper to explore the utility and potential 
of the CFOR-Hermite scheme.

The rest of this paper is organized as follows. In Section II, 
we give a brief introduction to the DSC algorithm and Hermite 
kernel. The CFOR scheme based on the Hermite kernel is also 
described in this section. The time marching schemes are 
discussed. Extensive numerical results are presented in Section 
III to demonstrate the use of the CFOR-Hermite and to explore 
its resolution limit. Concluding remarks are given in Section IV.

\section{Theory and algorithm}

This section present the spatiotemporal discretizations of
a flow equation. Fourier analysis is carried out for both 
the low-pass and high-pass filters generated from the Hermite 
kernel. Time marching schemes are described after the 
introduction to the discrete singular convolution (DSC) 
and conjugate filter oscillation reduction (CFOR) scheme.

\subsection{DSC-Hermite algorithm and CFOR scheme}

Discrete singular convolution (DSC) is an effective approach 
for the numerical realization of singular convolutions, which 
occurs commonly in science and engineering. The DSC algorithm 
has been widely used in computational science in recent years. 
In this section we will give a brief review on the DSC algorithm 
before introducing a new kernel, the Hermite kernel. For more 
details of the background and application of the DSC algorithm 
in partial difference equations, the reader is referred to 
Refs. \cite{weijcp,weijpa,weicmame,WPW}.

In the context of distribution theory, a singular 
convolution can be defined by
\begin{eqnarray}
F(t)=(T * \eta)(t) = \int_{-\infty}^{\infty} T(t-x) \eta (x)dx, 
\label{dsc1} 
\end{eqnarray}
where $T$ is a singular kernel and  $\eta (x)$ is an element of 
the space of test functions. Interesting examples include singular
kernels of Hilbert  type,  Abel type and delta type. The former 
two play  important roles in the theory of analytical functions, 
processing of analytical signals, theory of linear responses 
and Radon transform. Since delta type kernels are the 
key element in the theory of approximation and the numerical 
solution of differential equations, we focus on the singular 
kernels of delta type
\begin{eqnarray}
T(x)= \delta^{(q)}(x),  ~~~~ (q=0, 1, 2,\cdots),  
\label{deltakn}
\end{eqnarray}
where superscript $(q)$ denotes the $q$th-order ``derivative'' 
of the delta distribution, $\delta (x)$, with respect to $x$, 
which should be understood as generalized derivatives of 
distributions. When $q=0$, the kernel, $T(x)=\delta(x)$, is 
important for the interpolation of surfaces and curves, including
applications to the design of engineering structures. For hyperbolic
conservation laws and Euler systems, two special cases, $q=0$ and 
$q=1$ are involved, whereas for the full Naiver-Stokes equations, 
the case of $q=2$ will be also invoked.  Because of its singular
nature, the singular convolution of Eq. (\ref{dsc1}) cannot 
be directly used for numerical computations. In addition, the 
restriction to the test function is too strict for most practical 
applications. To avoid the difficulty of using singular 
expressions directly in numerical computations, we consider a
sequence of approximations of the form
\begin{eqnarray}
\lim_{\alpha\rightarrow\alpha_0, {\rm or }\beta\rightarrow\beta_0,
\cdots}\delta^{(q)}_{\alpha,\beta,\cdots}(x) 
= \delta^{(q)}(x),  ~~~~ (q=0, 1, 2,\cdots),  
\label{deltakn2}
\end{eqnarray}
where ${\alpha,\beta,\cdots}$ are a family of parameters which 
characterize the approximation and ${\alpha_0,\beta_0}$ are 
generalized limits. At least one of the parameters in Eq. 
(\ref{deltakn2}) has to take a limit.  One of the most commonly 
used DSC kernels is the regularized Shannon kernel (RSK)
\begin{eqnarray}
 \delta_{\sigma, \Delta}(x-x_k)= {\sin{\pi\over\Delta}(x-x_k)\over
{\pi\over\Delta} (x-x_k)}\exp\left({(x-x_k)^2\over2\sigma^2}\right),
\label{RSK}
\end{eqnarray}
where $\Delta$ is the grid spacing and $\sigma$ characterizes
width of the Gaussian envelop. The regularity of the kernel 
approximation can be improved by using a Gaussian regularizer. 
Obviously, for a given $\sigma\neq 0$, the limit of 
$\Delta\rightarrow 0$ reproduces the delta distribution.
Numerical approximation to a function $f$ (low-pass filtering) 
and its $q$th derivative (high-pass filtering) are realized by 
the following discrete convolution algorithm
\begin{eqnarray}
f^{(q)}(x) \approx \sum_{k=-W}^{W} 
     \delta^{(q)}_{\sigma,\Delta}(x-x_k)f(x_k),
              \quad (q=0,1,2,\cdots),
\end{eqnarray}              
where $\{ x_k\}_{k=-W}^W$ are a set of discrete grid points which 
are centered at $x$, and $2W+1$ is the computational bandwidth, or 
effective kernel support, which is usually smaller than the entire 
computational domain.

In this work, we consider the Hermite function expansion of 
the delta distribution\cite{Korevaar,Hoffman,weijpa}
\begin{eqnarray}\label{Hermite}
 \delta_{\sigma,\Delta}(x-x_k)
       =\f{1}{\sigma}\mbox{exp}\l(-\f{(x-x_k)^2}{2\sigma^2}\r)
        \sum_{k=0}^{n/2}\l( -\f{1}{4} \r)^k 
        \f{1}{\sqrt{2\pi} k!} H_{2k}\l(\f{x-x_k}{\sqrt{2}\sigma}\r),
\end{eqnarray}
where $H_{2k}(x)$ is the usual Hermite polynomial. The delta
distribution is recovered by taking the limit $n\rightarrow\infty$.
The Hermite function expansion of the Dirac delta distribution was 
proposed by Korevaar\cite{Korevaar} over forty years ago and was 
considered by Hoffman {\it et. al.} for wave propagations\cite{Hoffman}. 
Further discussion of its connections to wavelets and the DSC algorithm 
can be found in Ref. \cite{weijpa}(b). Just like other DSC kernels, 
the expression (\ref{Hermite}) is a low-pass filter, see FIG. 
\ref{compfig1}. Its approximation to the idea low-pass filter
is effectively controlled by the parameter $\sigma$.
The $q$th order derivative is readily given by the  analytical 
differentiation of expression (\ref{Hermite}) with the use of
recursion relations
\begin{eqnarray}
&&H_k^{(2)}(x)-2xH^{(1)}_{k}(x)+2kH_k(x)=0; \\
&&H^{(1)}_{k}(x)=2kH_{k-1}(x).
\end{eqnarray}
In numerical applications, optimal results are usually obtained if 
the Gaussian window size $\sigma$ varies as a function of the central 
frequency $\pi/\Delta$, such that $r=\sigma/\Delta$ is a parameter 
chosen in computations. In FIG. \ref{compfig1}, we plot the 
frequency responses of both the regularized Shannon kernel (RSK) 
and Hermite kernel (HK). It is seen that the effective frequency 
band of the RSK is slightly wider than that of the HK. The most 
distinctive difference of these two kernels occurs in their low-pass 
filters. Compared with the RSK, the HK low-pass filters decay fast 
at high wavenumbers and thus better approximate the ideal low-pass 
filter. As a result, they are expected to remove high-frequency errors 
much more efficiently. The ramification  of these differences in 
numerical computations is examined extensively in the next section.

The essence of the CFOR scheme is the following. Since the DSC 
high-pass filters have a limited  effective frequency band, 
short wavelength numerical errors will be produced when 
these high-pass filters are used for numerical differentiations. 
Conjugate low-pass filters are designed to remove such errors.
The resulting solution should be reliable over the effective 
frequency band, i.e., the low wavenumber region on which both 
high-pass and low-pass filters are very accurate. A total 
variation diminishing (TVD) switch is used to activate the low-pass 
filter whenever the total variation of the approximate solution 
exceeds a prescribed criteria. By adjusting the parameter $r$ in
the low-pass filter which determines the effective frequency band  
with a given DSC kernel, the oscillation near the shock can
be reduced efficiently. Since the DSC-Hermite kernel is nearly 
interpolative, the low-pass filter in the CFOR scheme 
is implemented via prediction ($u(x_i)\rightarrow u(x_{i+\f{1}{2}})$)
and restoration ($u(x_{i+\f{1}{2}}) \rightarrow u(x_i)$) \cite{weigu,guwei}.
The CFOR scheme with the RSK has been successfully applied to 
a variety of standard problems, such as inviscid Burgers' equation, 
Sod and Lax problems, shock-tube problems, flow past a forward
facing step and the shock/entropy wave interaction\cite{weigu,guwei}.

\subsection{Time marching schemes}

In this work, we  study four types of systems: 2D incompressible 
Euler system, 1D linear advective equation, 1D compressible Euler 
system and 2D compressible Euler system. For the incompressible Euler 
system, we adopt a third-order Runge-Kutta scheme, where a potential 
function is introduced to update the pressure in each stage and 
the bi-conjugate gradient method is used to solve the Poisson 
equation of the potential function for the pressure difference.
More detail of the scheme can be found Ref. \cite{WPW,wan01}.  For 
all the other systems, a standard forth-order Runge-Kutta scheme 
is employed for time integrations.

The choice of the time increment needs special care to assure that 
the solution is  not only stable but also time accurate. To 
ensure the iteration of the 2D incompressible Naiver-Stokes 
system to be stable under a third-order Runge-Kutta scheme,  
the following CFL condition should be satisfied \cite{EShu}
\begin{eqnarray}
\Delta t \l[ \underset{i,j}{\mbox{max}} \l(\f{|u_{ij}|}{\Delta x} 
      + \f{|v_{ij}|}{\Delta y} \r) 
      + \f{2}{\mbox{Re}} \l( \f{1}{\Delta x^2} 
      +\f{1}{\Delta y^2} \r)\r] \le 1.
     \label{dt1}
\end{eqnarray}
For the time integration of the linear advection equation, Hu 
{\it et al} \cite{hu94} demonstrated that a stricter restraint 
should be posed on $\Delta t$ to guarantee the solution to be 
time accurate. In this work, we test the numerical resolution 
of the proposed CFOR-Hermite scheme by using different small 
time increments $\Delta t$  to ensure the overall numerical 
error is dominated by the spatial discretization.

\section{Numerical Experiments}

In this section, a few benchmark numerical problems, including, 
the incompressible 2D Euler flow, a 1D wave equation, 1D and 2D 
compressible Euler systems, are employed to test shock-capturing 
capability of the proposed CFOR-Hermite scheme and to examine 
the numerical resolution of the scheme. We choose the CFOR-Hermite 
parameters of $n=88$ and $W=32$ for all the computations. 
The parameter of $r=3.05$ is used for high-pass filtering and
the prediction. For restoration, $r<3.05$ is used.

{\sc Example 1.} {\it Taylor problem}

As a benchmark test, the Taylor problem governed by the 2D 
incompressible Navier-Stokes equations is commonly used to access 
the accuracy and stability of schemes in CFD\cite{EShu,wan01}. 
The extension of this problem to the test of resolution can be 
found in Ref. \cite{yang01}. Since the viscous term leads to the 
dissipation of the initial wave, the long-time evolution accuracy 
of a numerical scheme can be masked by the viscous decay. 
In the present computation,  we consider the incompressible 
Euler version of this problem, i.e. Taylor problem without the 
viscous term.

\begin{figure}
 \centering  
 \subfigure[RSK]{
      \label{fig:subfig:a2}
      \includegraphics[width=8cm,height=8cm]{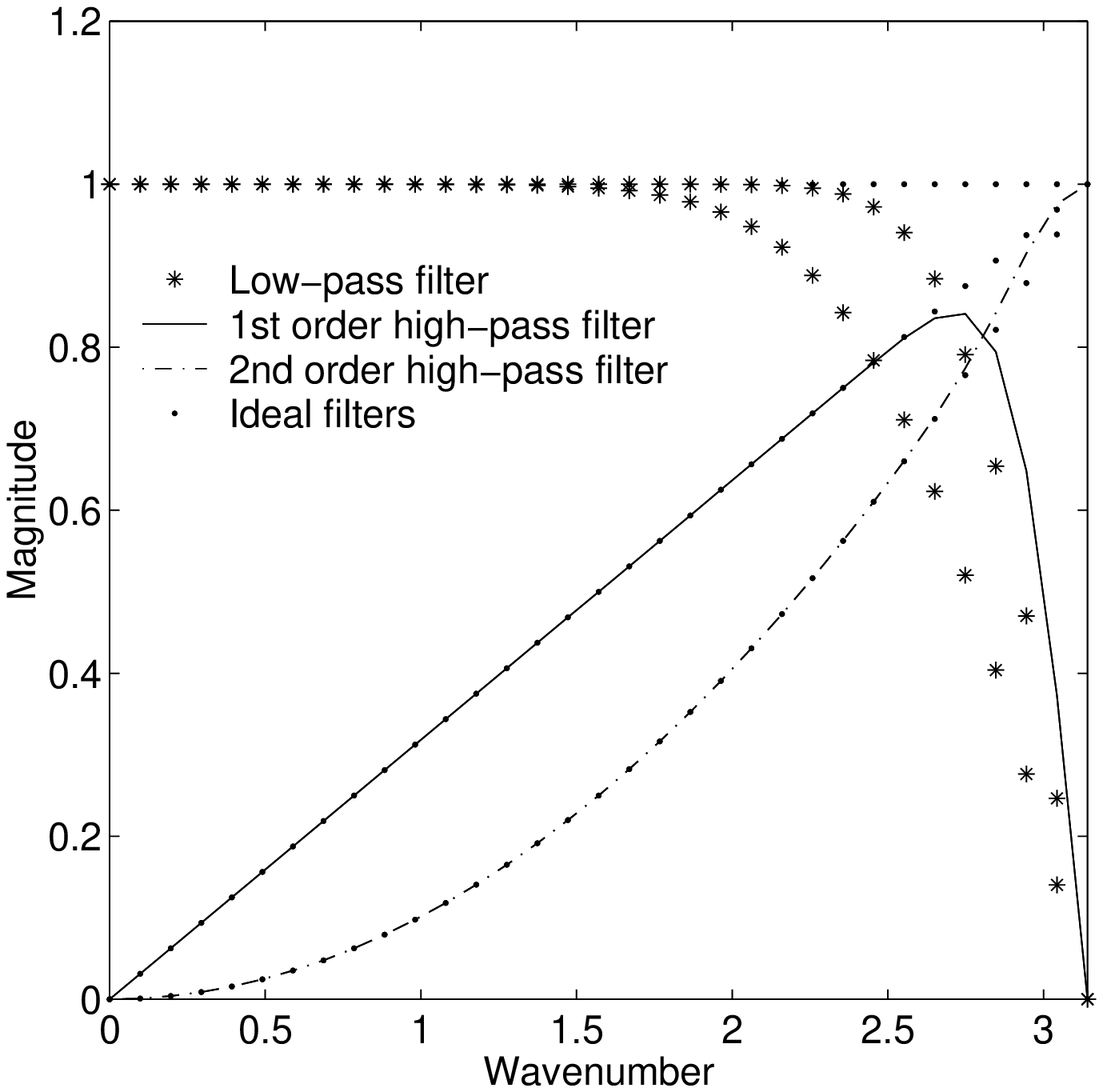}} 
 \subfigure[HK]{
      \label{fig:subfig:a1}
      \includegraphics[width=8cm,height=8cm]{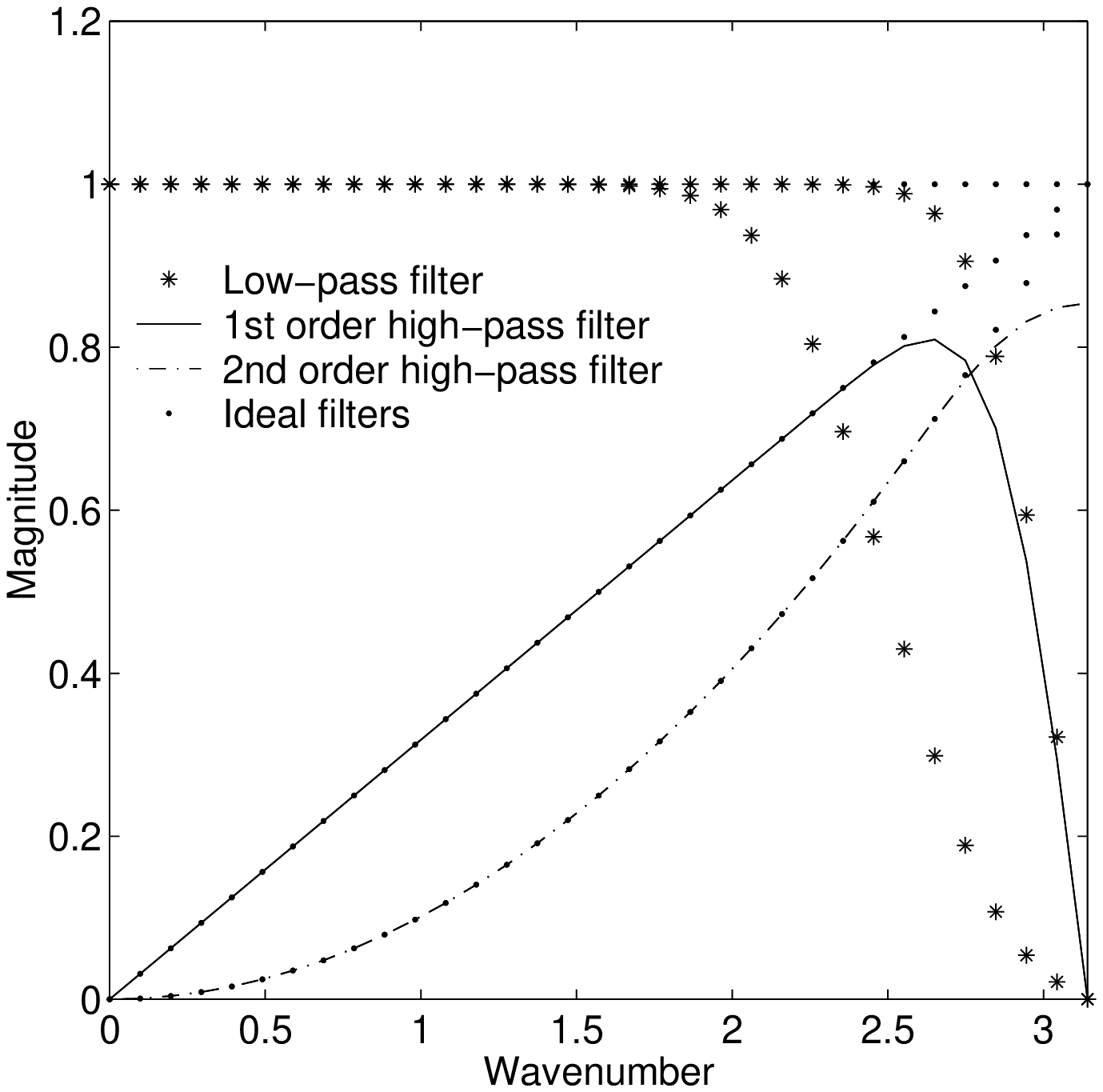}}
 \caption{Frequency responses of the conjugate DSC filters 
         (in the unit of $1/\Delta$),
          The maximum amplitude of the filters is normalized 
          to the unit. For HK: $W$=32; $r$=3.05
          for high-pass filters; $r$=2.5 for the left low-pass filter;
          $r$=3.05 for the right low-pass filter.
          For RSK: $W$=32; $r$=5.4
          for high-pass filters; $r$=2.0 for the left low-pass filter;
          $r$=3.2 for the right low-pass filter.} 
     \label{compfig1}
 
\newpage
 \subfigure[RSK]{
      \label{fig:subfig:II}
      \includegraphics[width=8cm,height=8cm]{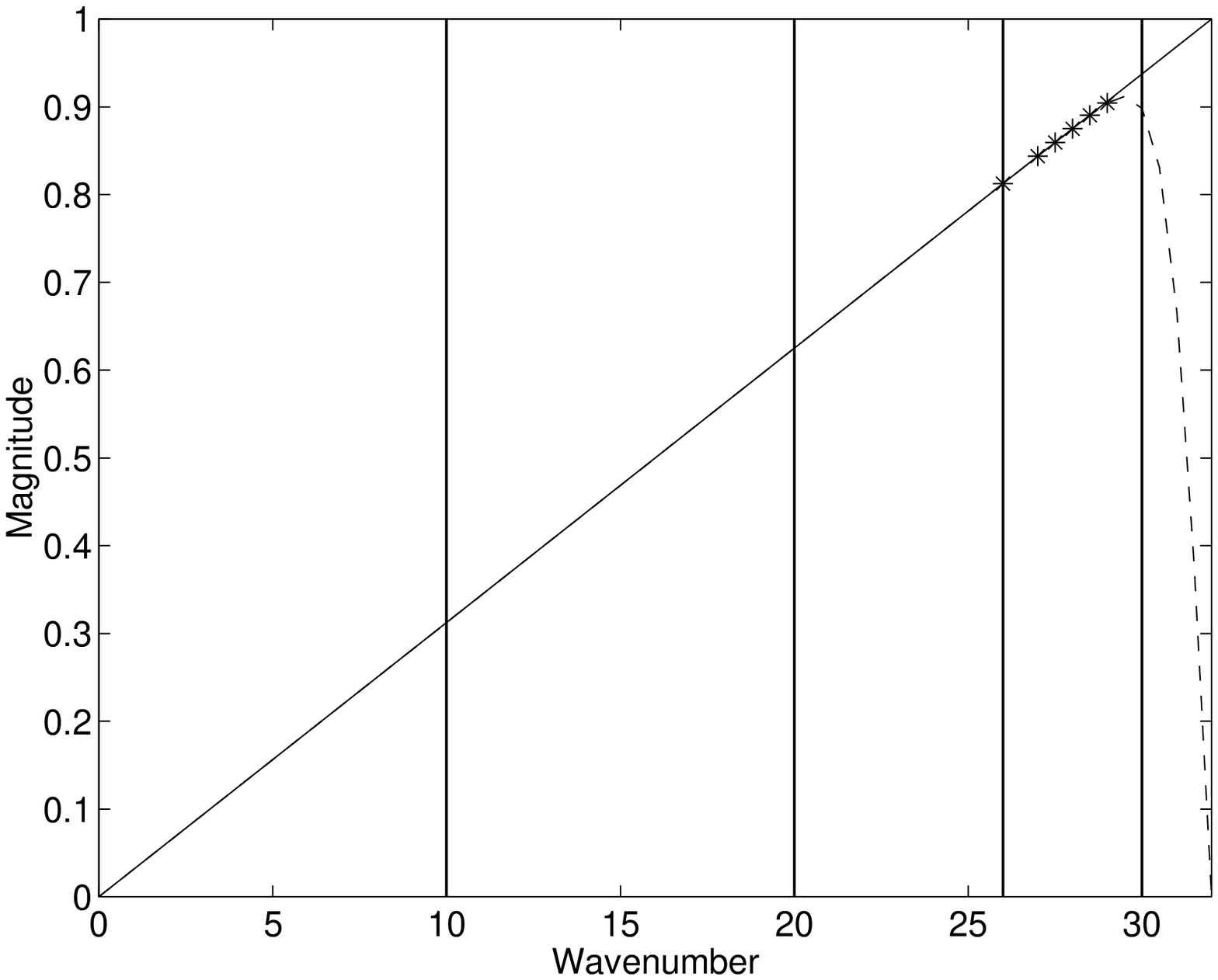}} 
 \subfigure[HK]{
      \label{fig:subfig:I}
      \includegraphics[width=8cm,height=8cm]{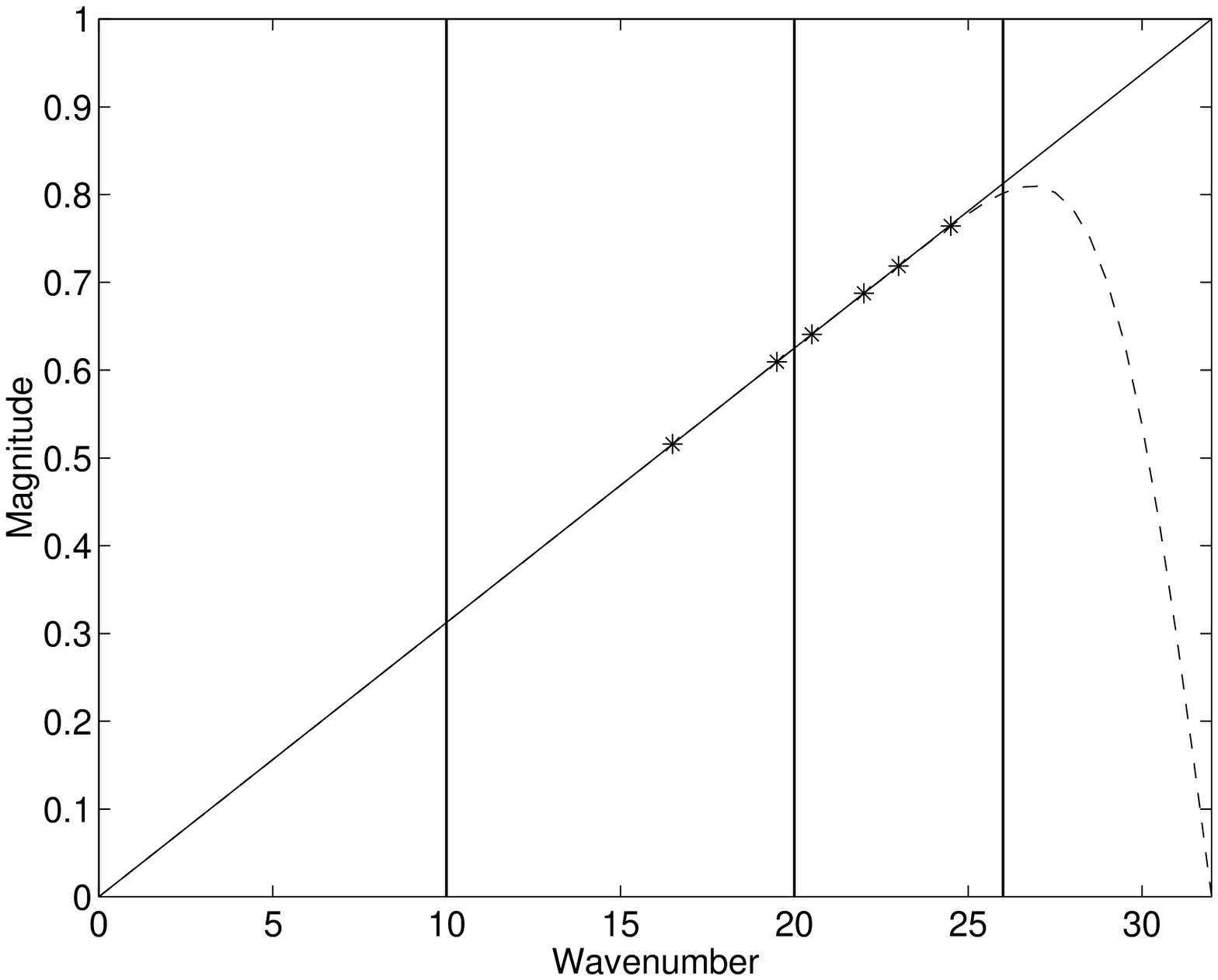}}
 \caption{Fourier analyses for the Taylor problem. Wavenumbers 
  are rescaled such that the maximum  is the half of mesh 
  number, i.e., 32 unit length. Sharp peaks are 1D projections of 
  frequency responses of the exact solutions.
  Six stars indicate the critical wavenumber
  beyond which  
  $\Delta\omega$, will be greater than 
  $10^{-10},10^{-7},10^{-6},10^{-5},10^{-4},10^{-3}$, respectively
  (from the left to the right). 
  Here, $\Delta\omega$ is the difference between the exact 
  first-order derivative (the diagonal line) 
  and the kernel approximation (the   dashed line).}
  \label{taylorfig1} 
 
 \newpage

 \includegraphics[width=16cm]{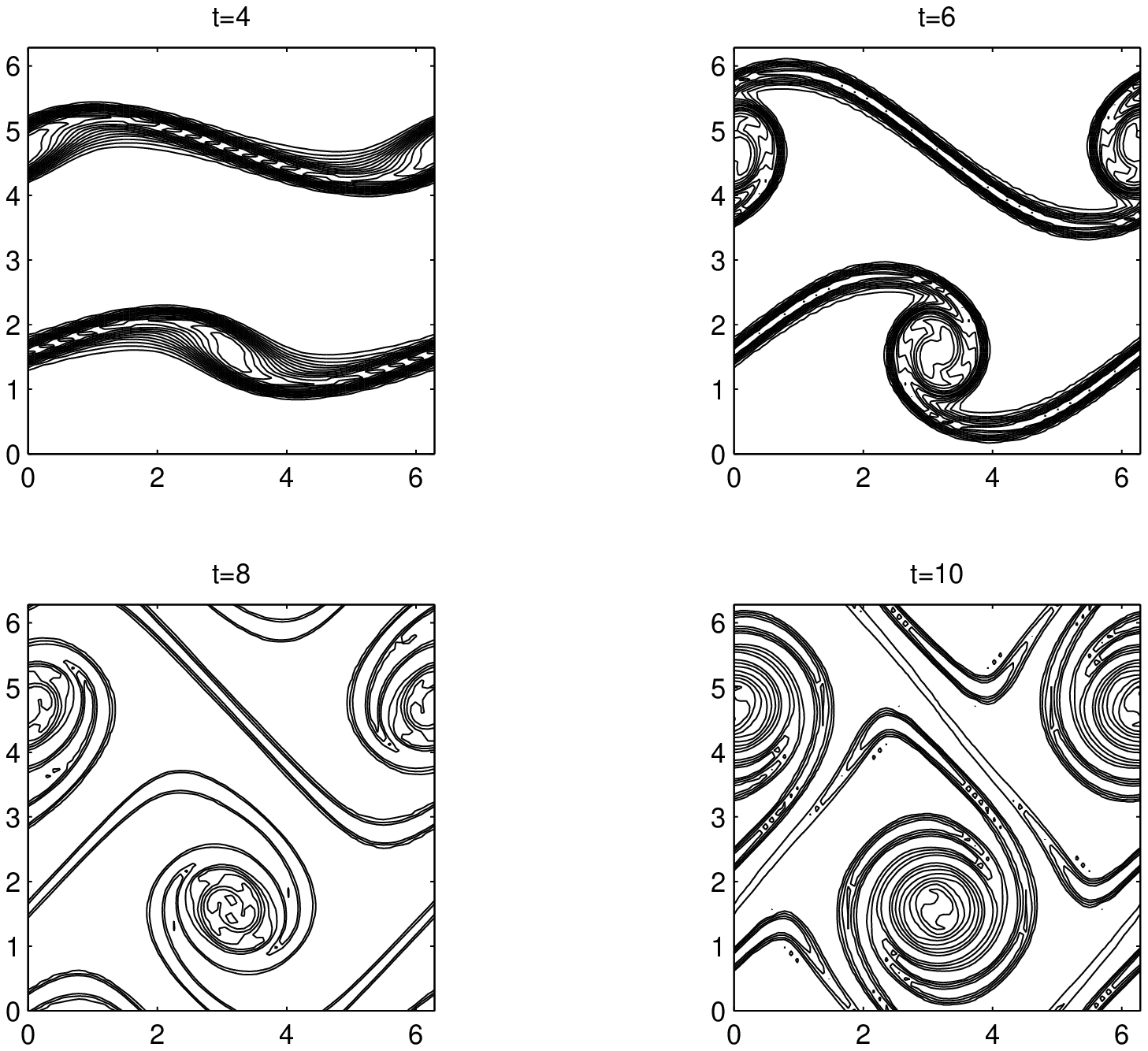}
 \caption{Contours of vorticity for double shear layer problem at 
          different time.} 
 \label{shearfig1}

 \newpage

 \includegraphics[width=12cm]{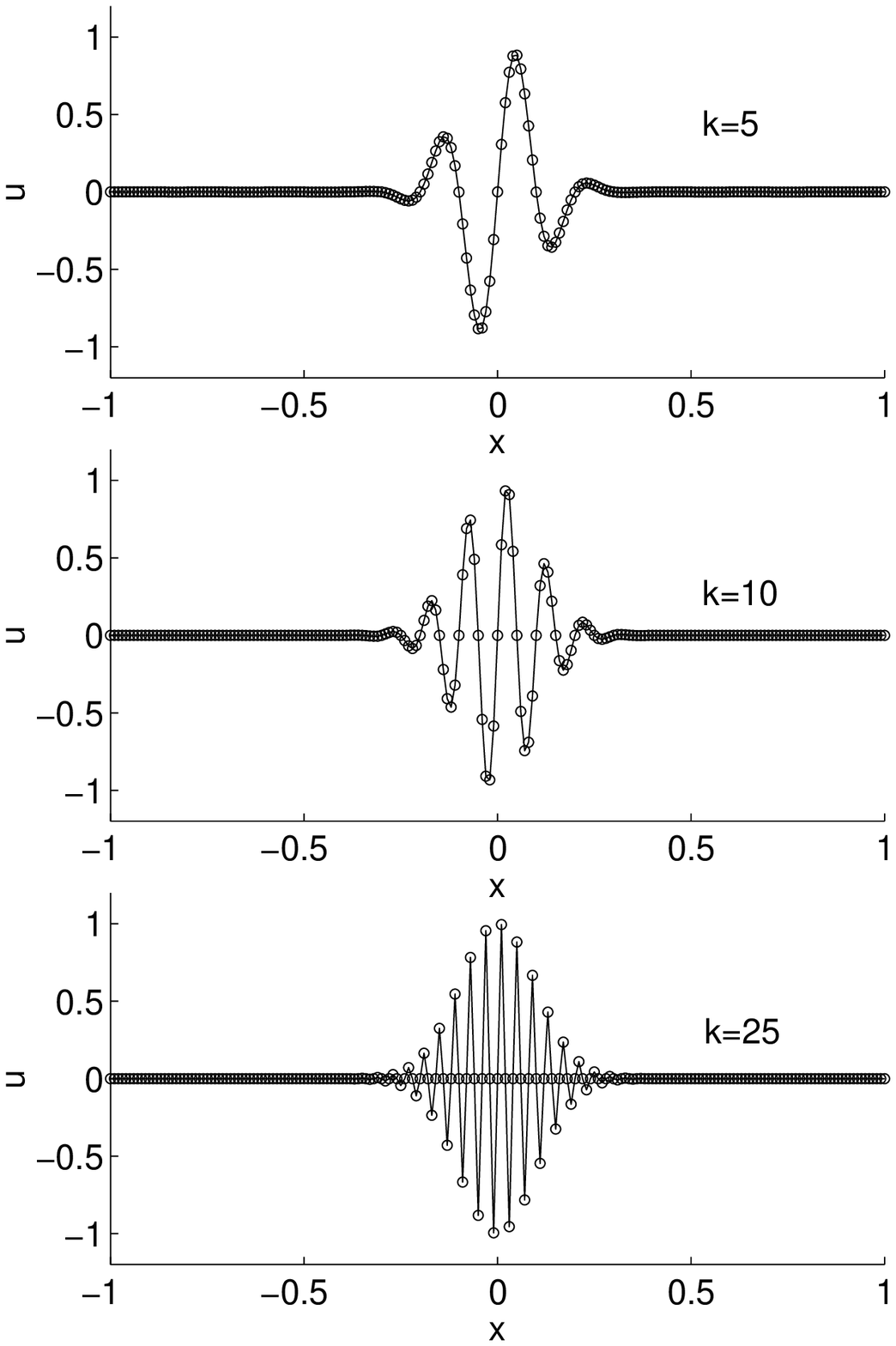}
 \caption{The sine-Gaussian wavepackets at $t=100$. 
  Exact solutions are plotted in solid lines and numerical
   ones are given by circles.
    }
  \label{wavepacket}

 \newpage

 \includegraphics[width=12cm]{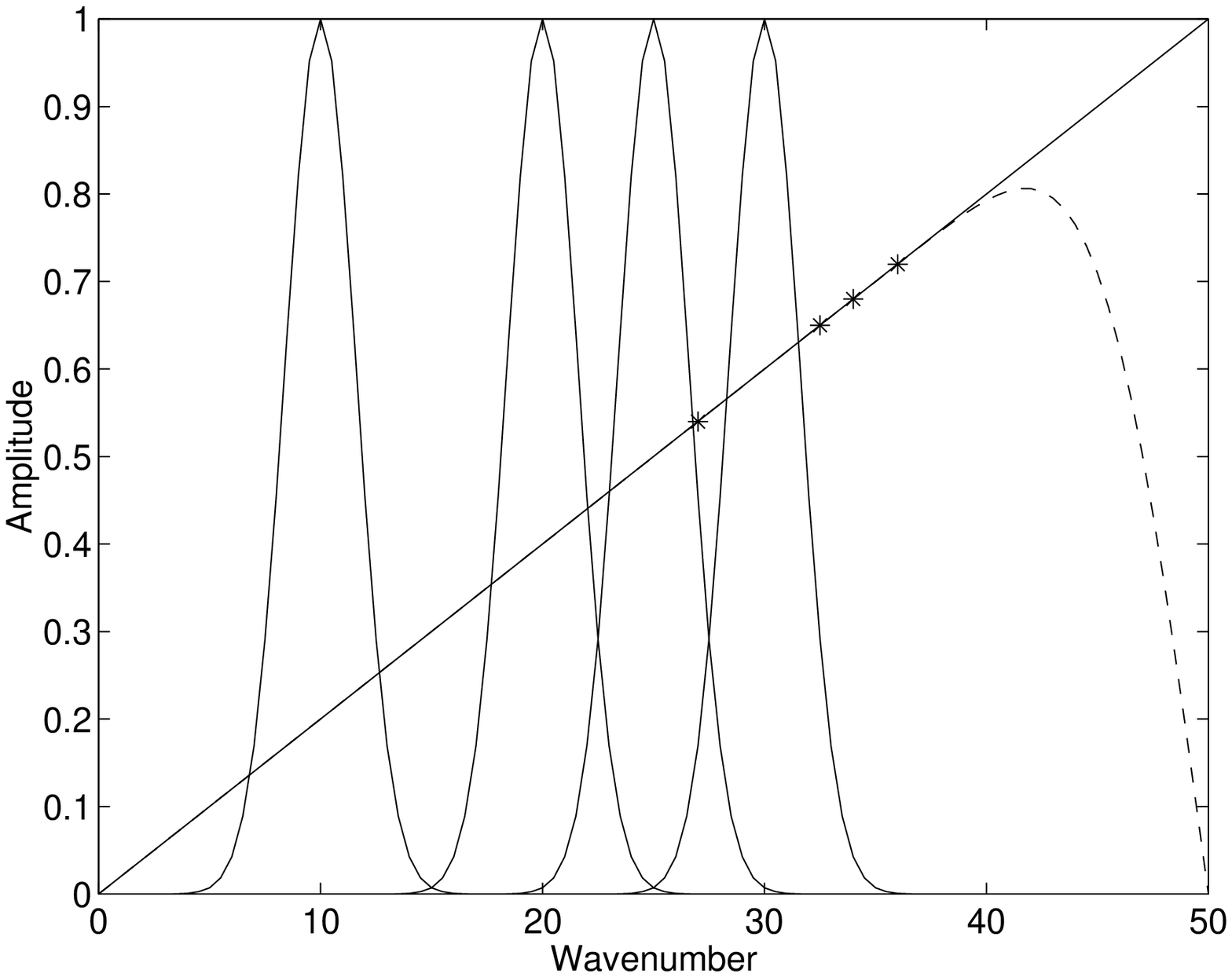}
 \caption{Fourier analyses for the advective sine-Gaussian wavepacket.
  Wavenumbers are rescaled such that the maximum  is the half of mesh 
  number, i.e., 50 unit length.  Four Gaussian peaks correspond to 
  $k=10,20,25,30$,  respectively. Four stars indicate the critical 
  wavenumbers beyond which 
  $\Delta\omega$ will be greater than
  $10^{-10},10^{-6},10^{-5},10^{-4}$, respectively 
  (from the left to the right). 
  Here, $\Delta\omega$ is the difference between the exact 
  first-order derivative (the diagonal line) 
  and the kernel approximation (the   dashed line).} 
  \label{advcfig1}
 
\newpage

 \subfigure[$t=2$]{
     \label{fig:subfig:A}
     \includegraphics[width=8cm]{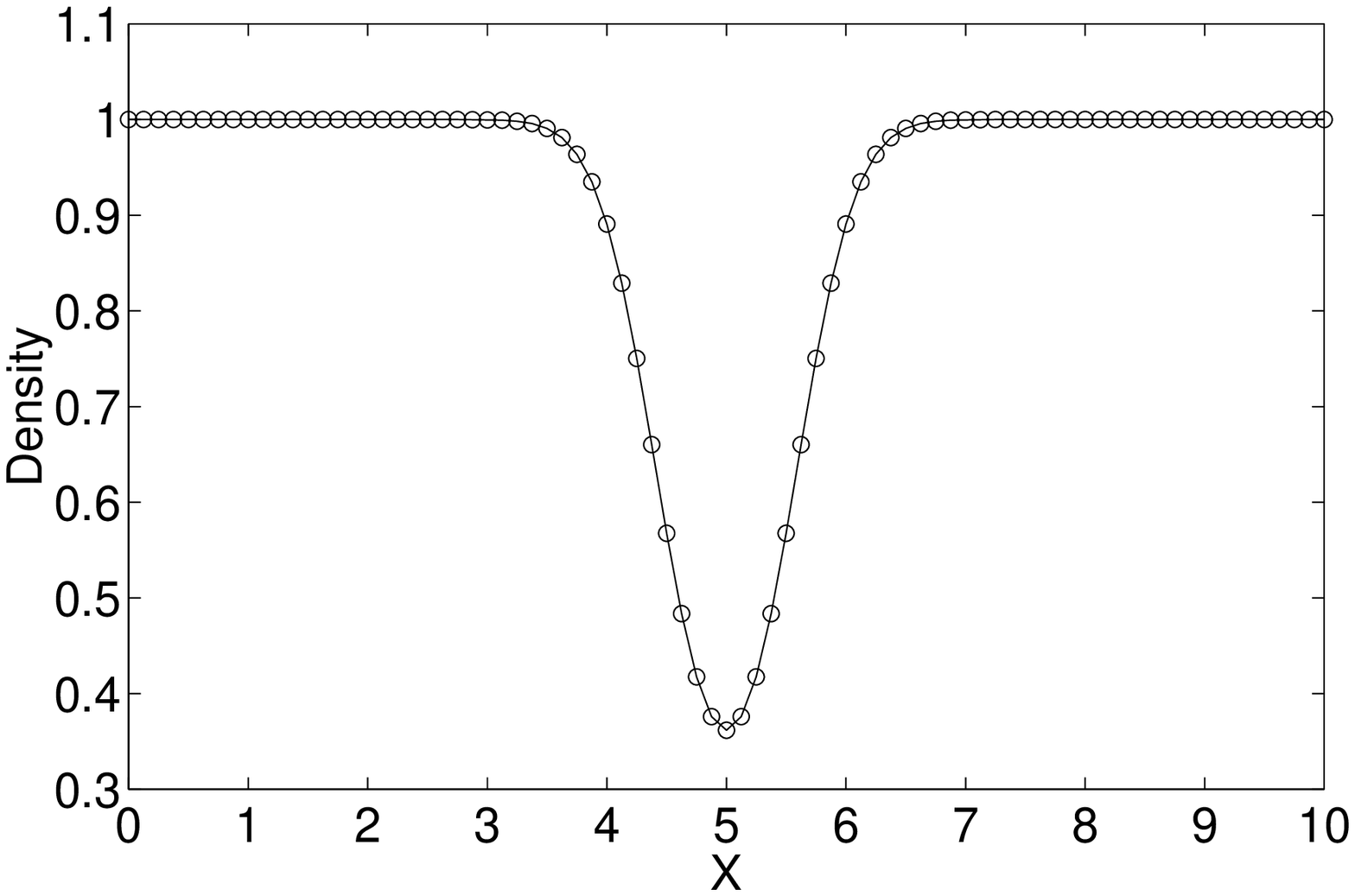}}
 \hspace{0.5cm} 
 \subfigure[$t=10$]{
     \label{fig:subfig:B}
     \includegraphics[width=8cm]{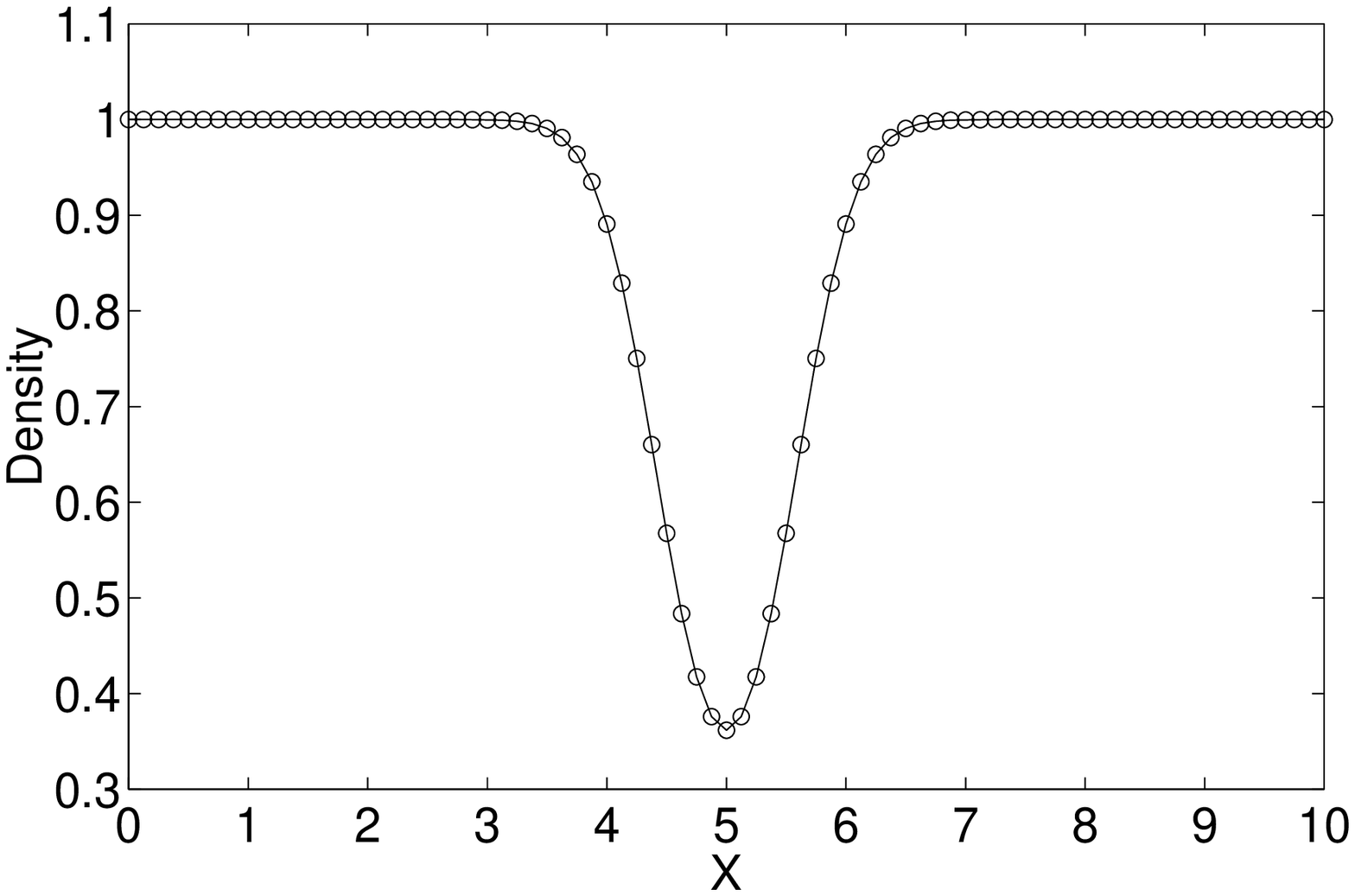}}  \\ 
  \subfigure[$t=50$]{
     \label{fig:subfig:C}
     \includegraphics[width=8cm]{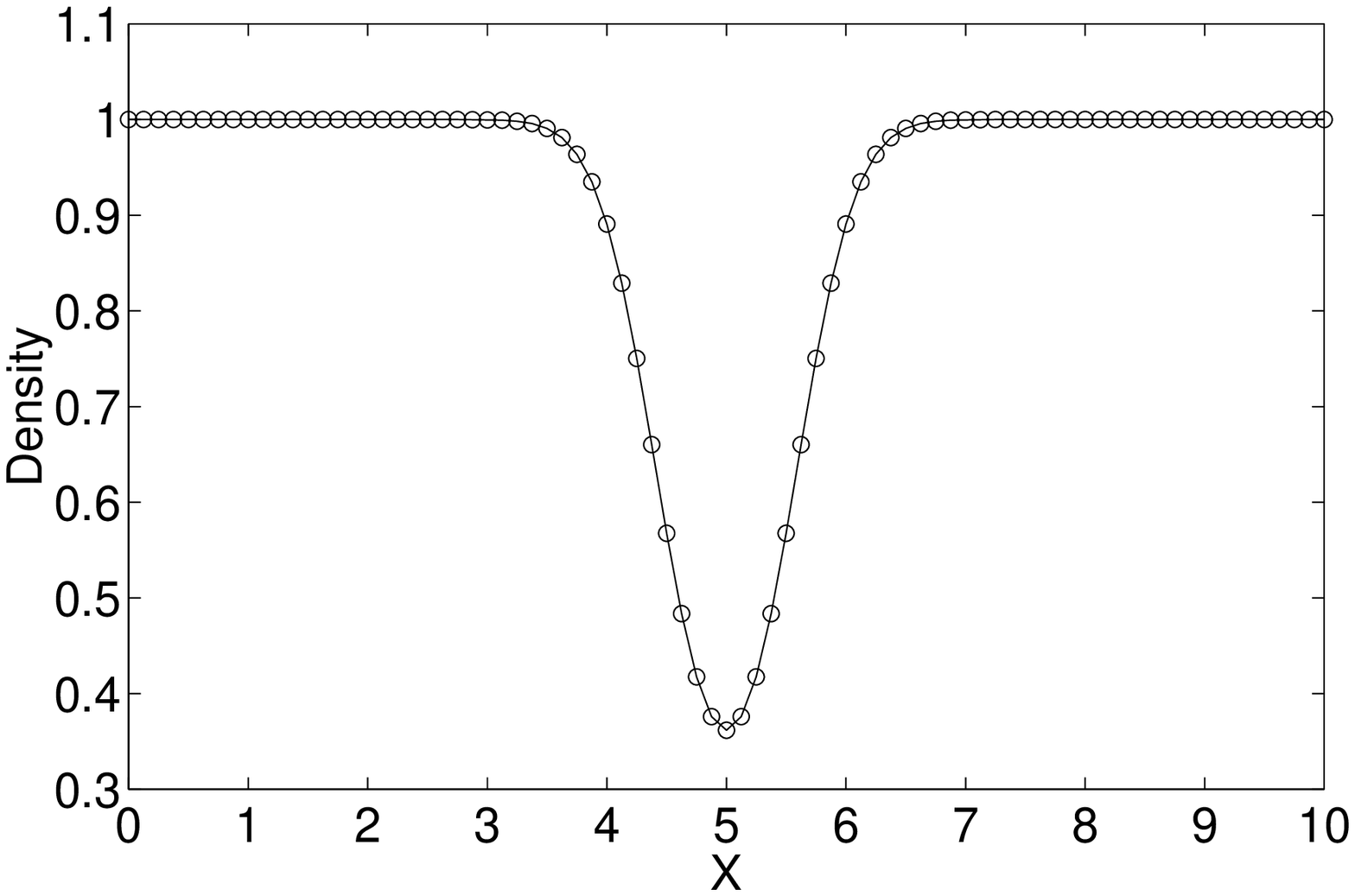}}
 \hspace{0.5cm}
 \subfigure[$t=100$]{
     \label{fig:subfig:D}
     \includegraphics[width=8cm]{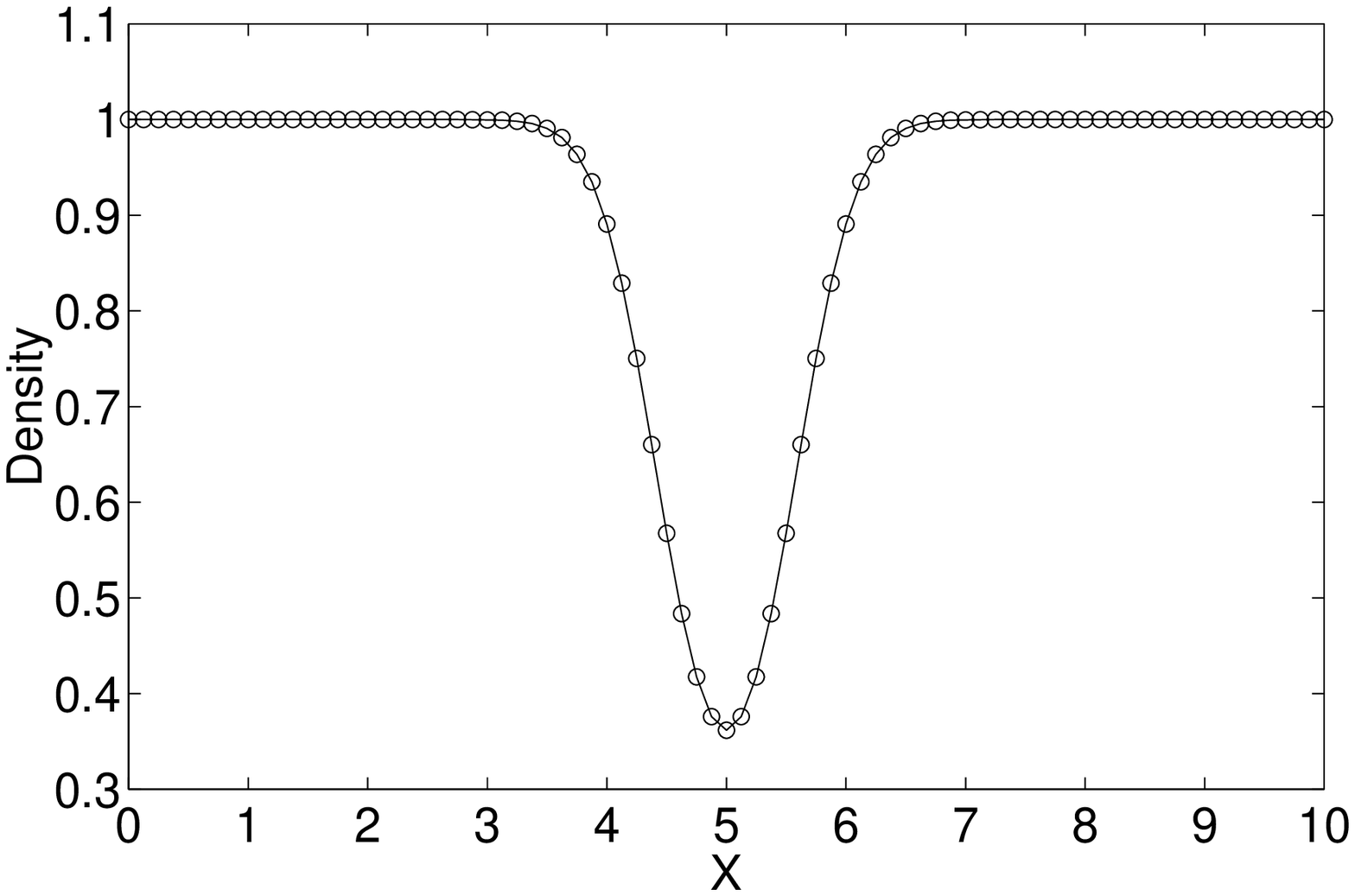}}
 \caption{Density profiles in horizontal cutting at four times. 
          Solid line is the exact profile and the circle dots 
          denote the numerical results. $N=80$; CFL=0.5. 
          See also the Table \ref{vortextable3} for numerical errors.} 
          \label{vortexfig1}
 
 \subfigure[$k=13, N=400$]{
     \label{fig:subfig:1}
     \includegraphics[width=12cm]{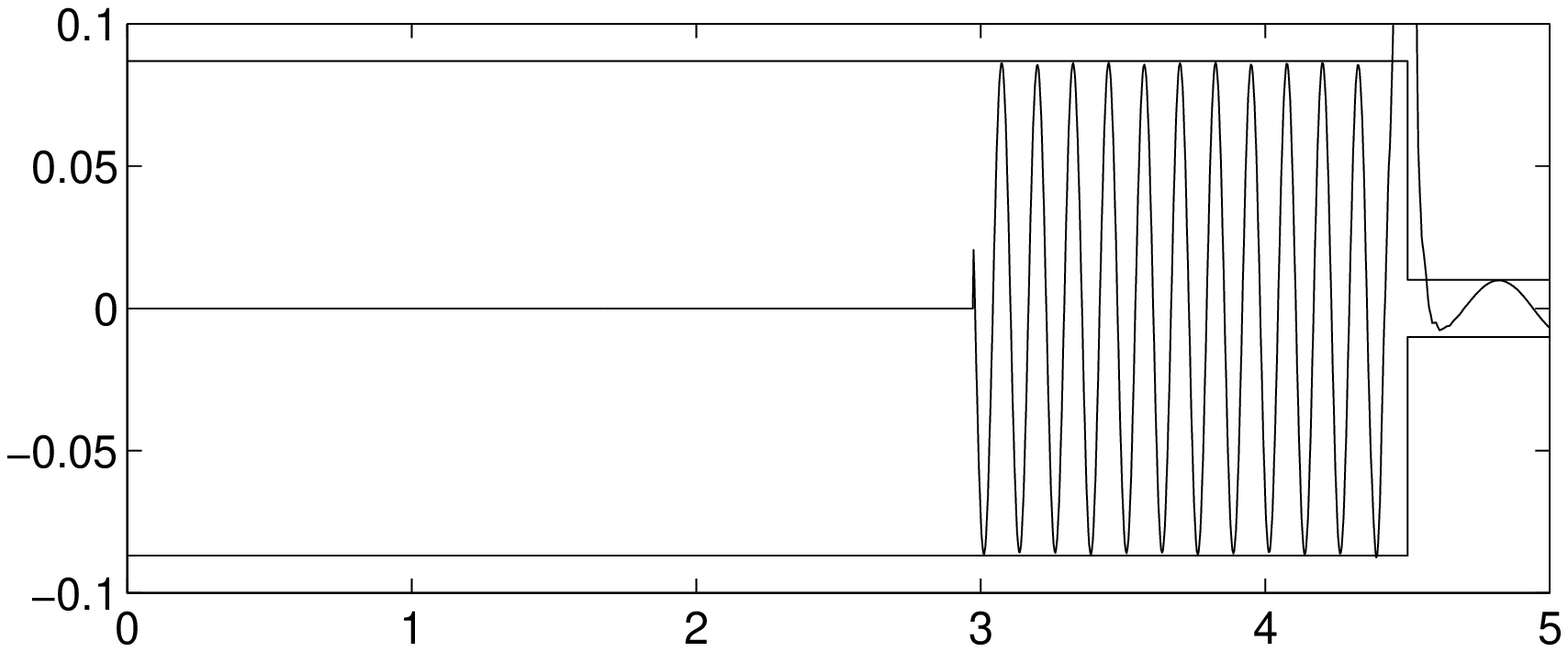}} 
 \subfigure[$k=13, N=800$]{
     \label{fig:subfig:2}
     \includegraphics[width=12cm]{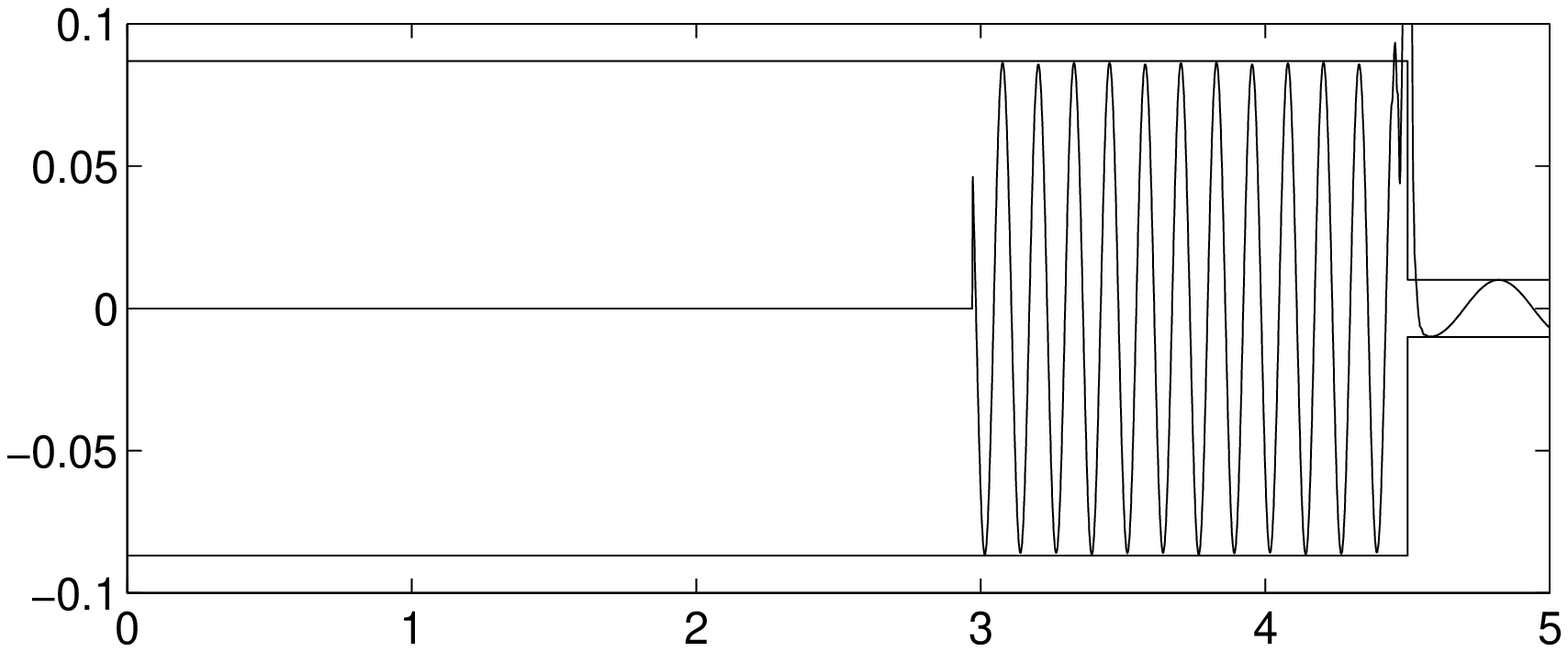}}
 \subfigure[$k=26, N=400$]{
     \label{fig:subfig:3}
     \includegraphics[width=12cm]{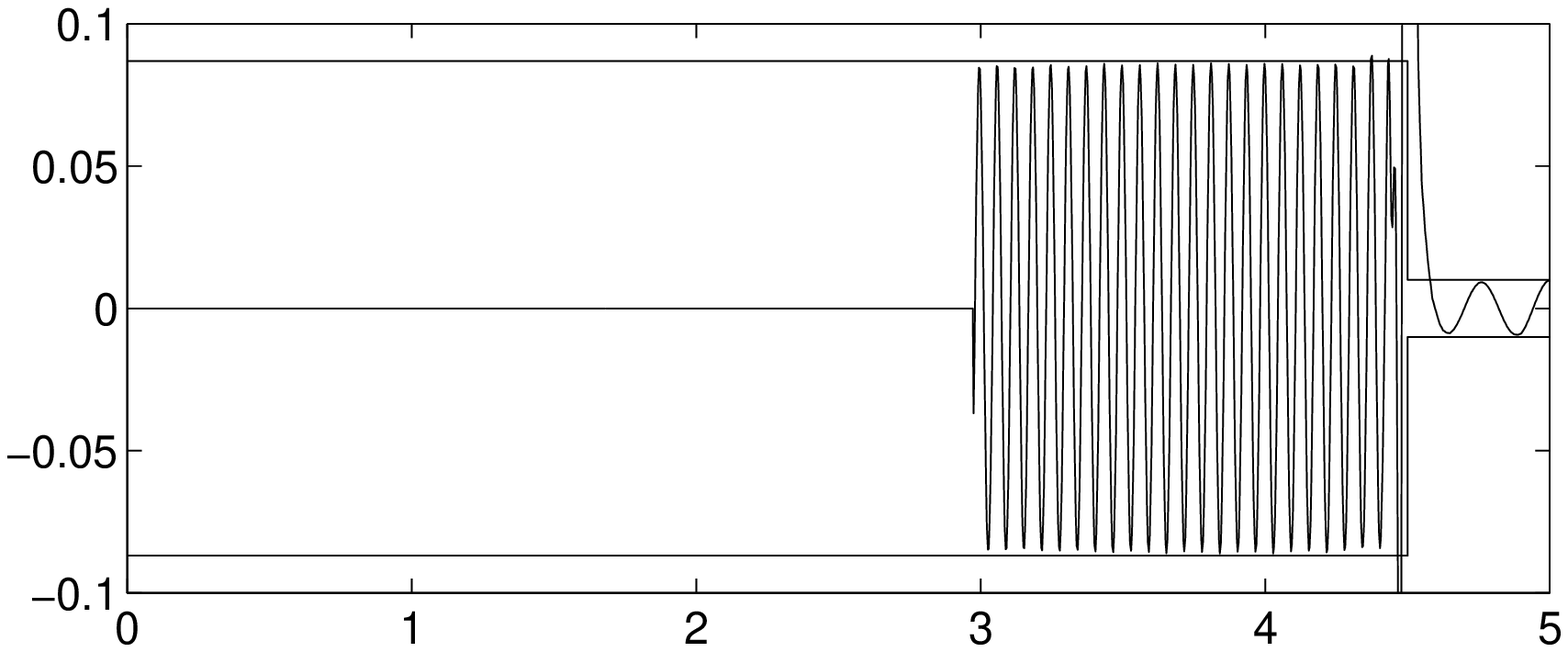}}
 \newpage 
 \subfigure[$k=26, N=800$]{
     \label{fig:subfig:4}
     \includegraphics[width=12cm]{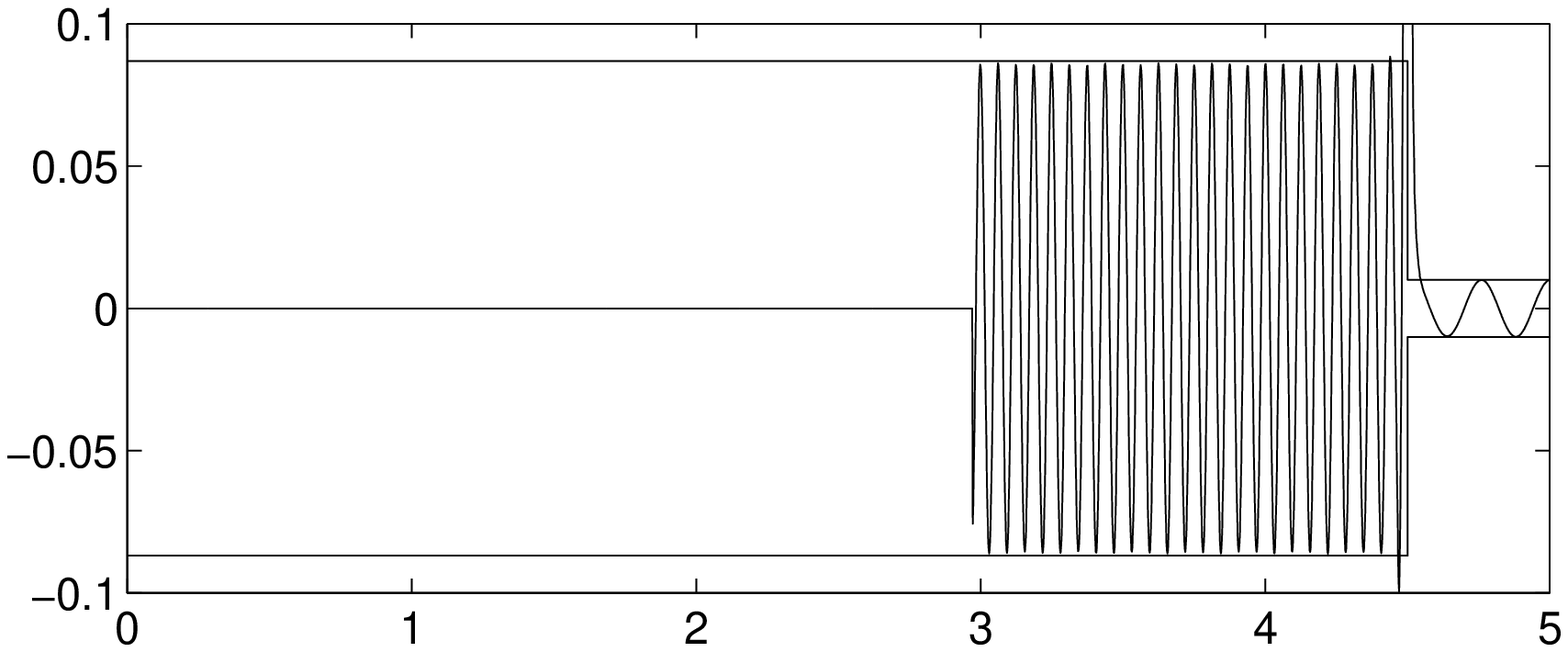}}
 \subfigure[$k=52, N=800$]{
     \label{fig:subfig:5}
     \includegraphics[width=12cm]{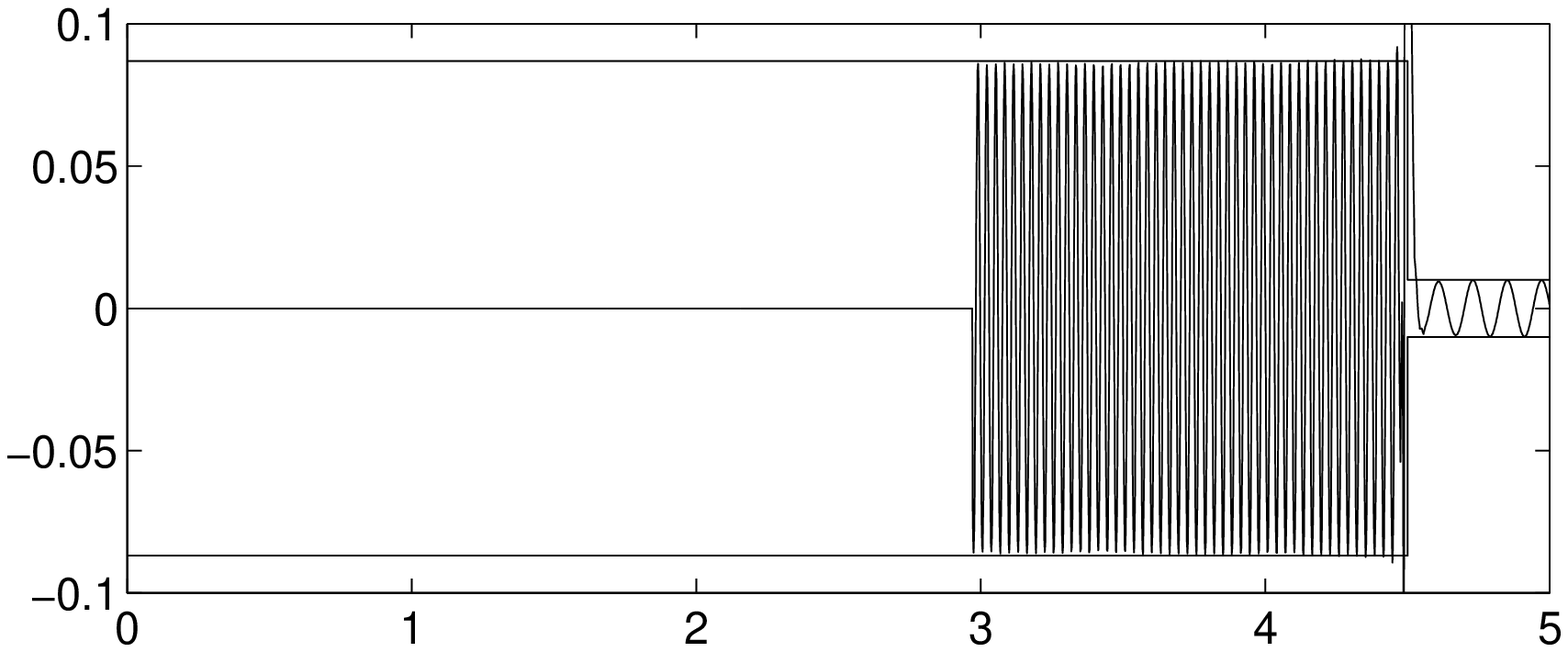}} 
 \subfigure[$k=52, N=1200$]{
     \label{fig:subfig:6}
     \includegraphics[width=12cm]{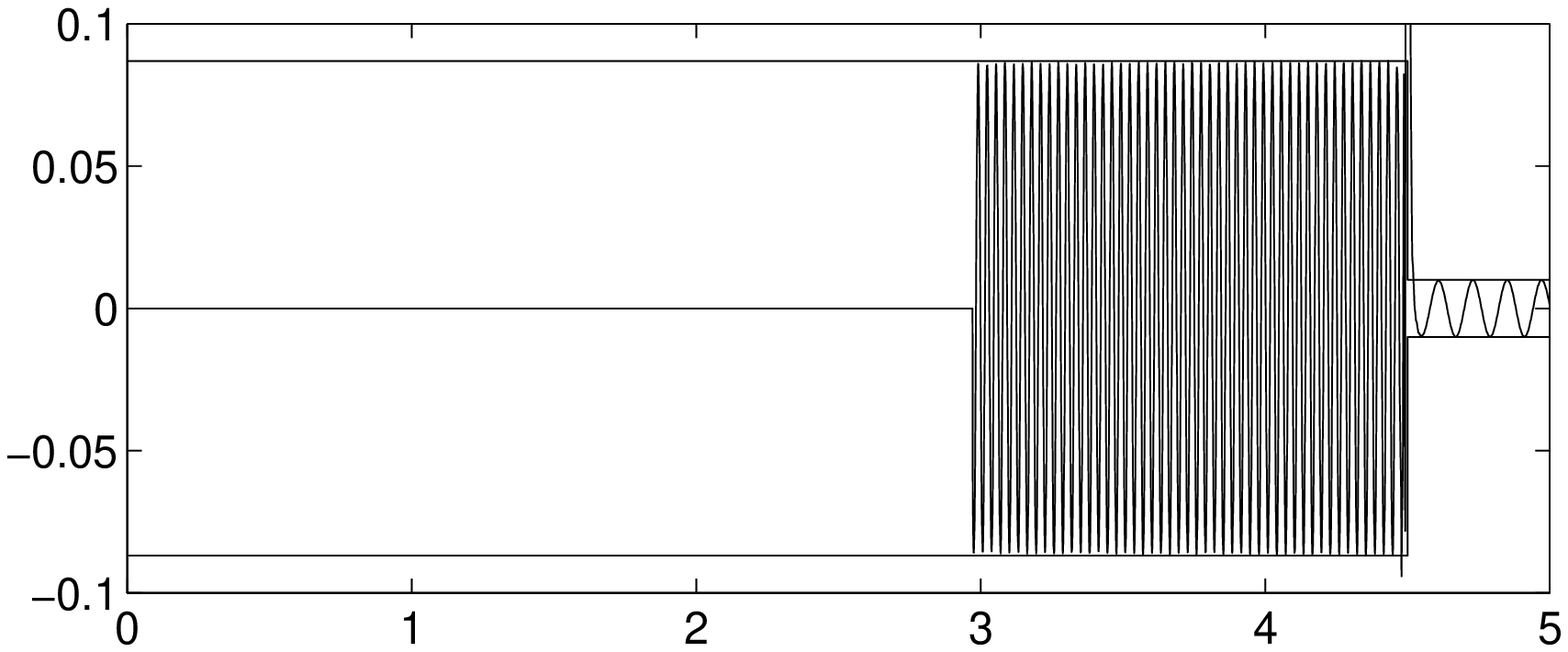}}
 \newpage
 \subfigure[$k=65, N=1000$]{
     \label{fig:subfig:7}
     \includegraphics[width=12cm]{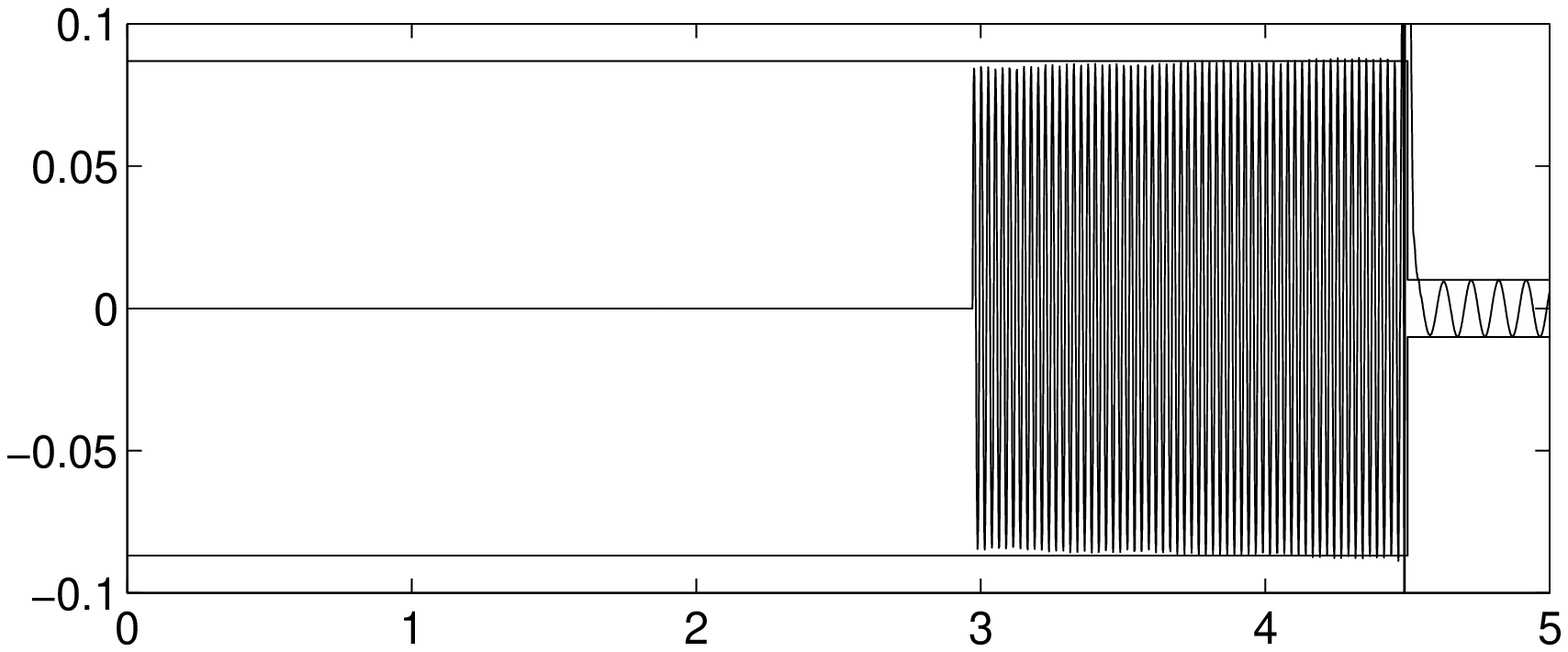}} 
 \subfigure[$k=65, N=1200$]{
     \label{fig:subfig:8}
     \includegraphics[width=12cm]{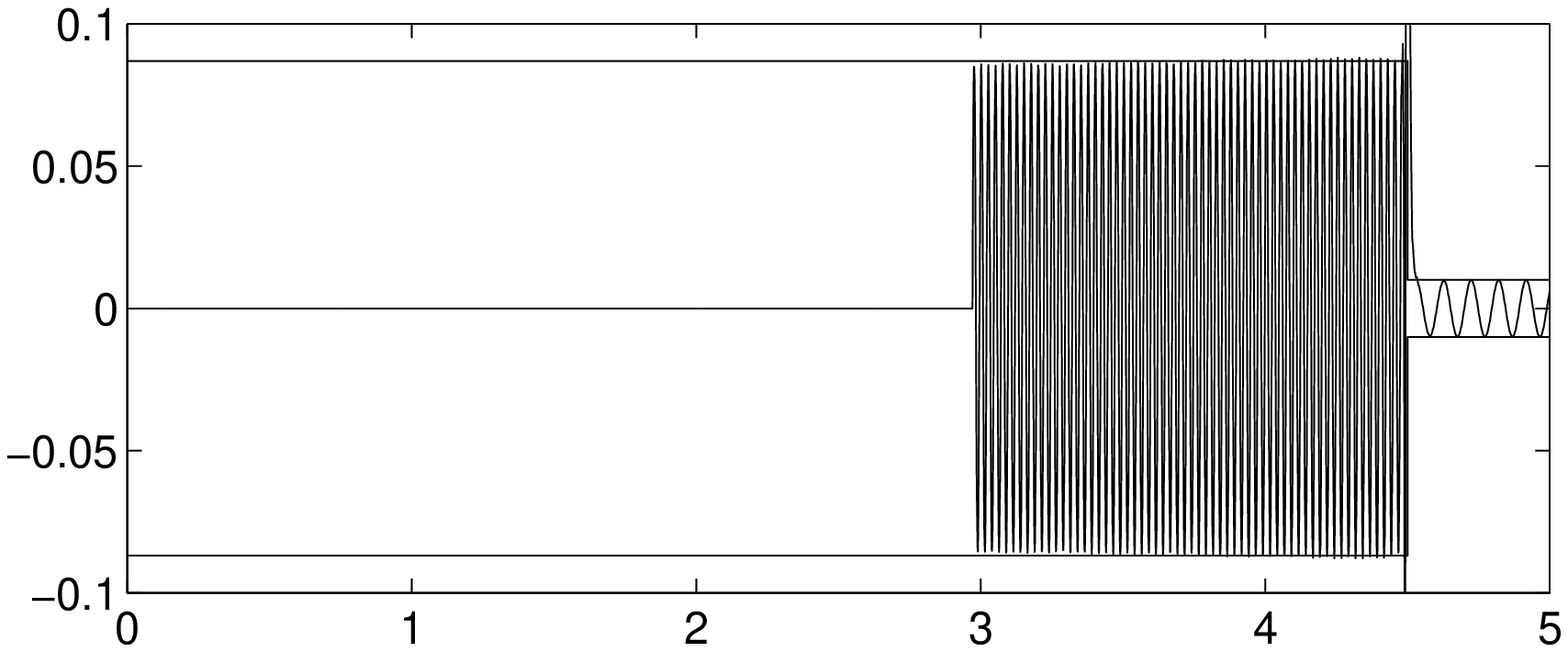}}
 \subfigure[$k=70, N=1200$]{
     \label{fig:subfig:9}
     \includegraphics[width=12cm]{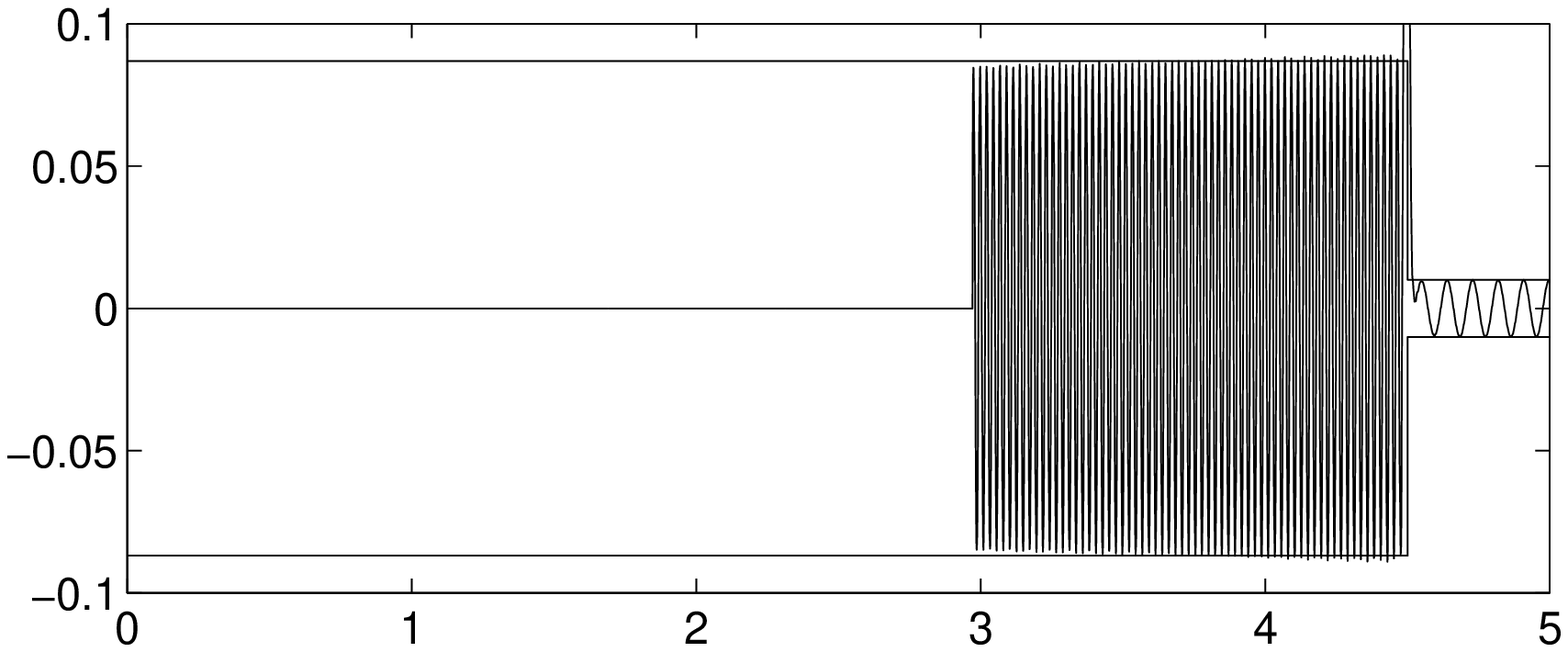}}
 \caption{Entropy waves for different pre-shock wavenumber.} 
  \label{shockfig1} 
\end{figure}

With a periodic boundary condition, the 2D time dependent 
incompressible Euler equation admits the following exact solution
\begin{eqnarray}
u(x,y,t)  & = &  -\mbox{cos}(kx)\mbox{sin}(ky)     \\
v(x,y,t)  & = &  \quad \mbox{sin}(kx)\mbox{cos}(ky)  \\
p(x,y,t)  & = &  -\l( \mbox{cos}(2kx)+\mbox{cos}(2ky) \r)/4.0,
\end{eqnarray} \label{taylorsol}
where $u,v,p$ denote the velocity in $x$- and $y$-directions and 
the pressure, respectively. Here $k$ is the wavenumber of the 
velocity solution. Note that the wavenumber of the pressure is 
twice as much as that of the velocity. Therefore, for a given $k$ 
value, the overall accuracy is limited by the capability of resolving 
the pressure. The computation domain is chosen as $[0,2 \pi] \times 
[0, 2\pi]$ with periodic boundary conditions in both directions. 
Mesh size of the computational domain is $1/N$, and $N=64$ is 
used in all the tests in this problem. We examine the resolution 
of the CFOR scheme with a variety of $k$ values. The numerical 
errors with respect to the exact velocity $u$ are listed in 
Table \ref{taylortab1}, whereas the ratio of grid points per 
wavelength (PPW) is computed according to the pressure wavenumber 
($2k$).  The PPW is a good measure of the numerical resolution 
for many problems of wave nature. Comparison is also carried out 
between the regularized Shannon kernel (RSK) and the Hermite kernel, 
since we speculate that for this smooth problem the RSK should 
have a better performance with a high wavenumber due to its wider 
effective frequency band as shown in FIG. \ref{taylorfig1}.

The results in Table \ref{taylortab1} exhibit high accuracy of 
the DSC algorithm, realized both by the RSK and HK. It is not 
surprising that with the increase of the wavenumber $k$, the
accuracy of the numerical scheme degrades due to the finite
effective frequency band. In FIG. \ref{taylorfig1}, we analyze  
the Fourier resolution of the RSK and HK (the stars) against 
the Fourier response of the exact solution  with a given $k$ 
(the straight lines). Their relative locations in the Fourier 
domain tell us the approximation errors. The accuracy of the 
final results is eventually limited by the Fourier resolution, 
i.e., the effective frequency band of the kernel. First, we note 
that results in Table \ref{taylortab1} are qualitatively 
consistent with the analysis in FIG. \ref{taylorfig1}.
For example, the analysis indicates that the RSK has a better
Fourier resolution than the RK. The numerical results in Table 
\ref{taylortab1} show the same tendency.  Quantitatively, 
the frequency response error of the first order derivative 
approximated by using the RSK is about $10^{-10}$ at wavenumber 
$2k=26$. The errors of numerical results are  also about 
$5\times10^{-11}$ and thus match the Fourier analysis. In general, 
there is a good consistence between the prediction of Fourier 
analysis and the experimental results for the RSK. In a dramatic 
case, the PPW ratio is reduced to 2.1, which is very close to 
the Nyquist limit of 2. It is noted that the CFOR-RSK still 
performs very well for this case. This confirms the spectral 
nature of the CFOR-RSK scheme.

For the Hermite kernel, the errors at $2k=20$ and 26 are slightly 
smaller than those predicted by the Fourier analysis. This 
discrepancy is due to the fact that the actual errors measured 
are for the velocity rather than for the pressure. The pressure 
error might not {\it fully} propagate to the velocity because 
only the pressure gradient contributes to velocity changes.
Comparing with the RSK, the Hermite kernel has a narrow effective
frequency band, and thus performs not as well as the RSK. 
However, a careful examination on the FIG. \ref{compfig1} 
reveals that the Hermite kernel has a sharp decay in its
conjugate low-pass filter as the wavenumber increases. Therefore,
it is expected to perform better in the case of shock-capturing 
involving natural high frequency oscillations. This is indeed
the case as demonstrated in the example of shock/entropy wave 
interaction.

{\sc Example 2.} {\it Double shear layers}

This problem is also governed by the 2D incompressible Euler 
equations with periodic conditions\cite{EShu}. The problem is 
not analytically solvable due to the following initial conditions
\begin{eqnarray}
u(x,y,t) & = &  \begin{cases}
   \mbox{tanh} \l(\f{2y- \pi}{2 \rho} \r), & \mbox{if} ~y \le \pi \\
   \mbox{tanh} \l(\f{3\pi-2y}{2 \rho} \r), & \mbox{if} ~y > \pi
                \end{cases} \\
v(x,y,t) & = & \delta \mbox{sin}(x),
\end{eqnarray} \label{shearsol}
where, $\rho$ represents the thickness of the shear layer. Such 
initial conditions describe a thin horizontal shear layer disturbed 
by a small vertical velocity. The initial flow quickly evolves 
into roll-ups with smaller and smaller scales until the full 
resolution is lost\cite{EShu}. A small $\rho$ speeds up the 
formation of small-scale  roll-ups. The initial shear layers
become discontinuous as $\rho$ approaches zero. In order to resolve 
the fine structures of roll-ups, high-order schemes are necessary. 
However, the use of high-order schemes along is not sufficient for 
achieving high resolution. The shock-capturing ability is required
for solving this problem. Bell, Colella and Glaz\cite{bell89} 
introduced  a second order Godunov-type shock-capturing scheme 
to simulate this incompressible flow. E and Shu\cite{EShu} 
considered an ENO scheme for handling this case. Recently,  Liu 
and Shu\cite{liu00} applied a  discontinuous Galerkin method to 
this problem. In the present work, we examine the performance of 
the CFOR-Hermite scheme for this problem. The low-pass filter will 
be activated whenever spurious oscillations are generated from 
the flow fields. The parameter $r=2.6$ is used in the restoration.

We consider the case of $\rho=1/15$ and $\delta=0.05$ in our
simulation. The time increment is fixed as $0.002$ and a mesh of 
$64^2$ is used. We note that a finer mesh was considered in previous
work\cite{bell89,EShu,liu00}. One of our objectives is to 
demonstrate that the CFOR gives excellent results with a coarse
mesh.

FIG. \ref{shearfig1} depicts the velocity contours at $t=4, 6, 
8$ and 10. It is seen that very good resolution is obtained in 
such a coarse grid. The present results are very smooth and the 
integration is very stable. Spurious oscillations are effectively 
removed by the conjugate low-pass filter. There is no serious 
distortion occurring to the channels connecting  individual
vorticity centers. It is worthwhile to stress that if a numerical 
scheme is highly dissipative, the vorticity cores which should be 
skew circle shape (see FIG. \ref{shearfig1} for cases  $t=4,6$) 
will be smeared to a circle shape. The quality of the present solution 
at early times is comparable to that of the spectral solution given 
by E and Shu\cite{EShu}.

{\sc Example 3.} {\it Advective sine-Gaussian wavepacket}

This 1D problem is designed to examine the resolution of 
the CFOR-Hermite approach with a linear advection equation
$\f{\p u}{\p t} + c \f{\p u}{\p x}=0$ with a smooth solution. 
In this problem, the initial data is taken as a wavepacket 
constructed by using a Gaussian envelop modulated by a
sine wave with a wavenumber $k$
\begin{eqnarray}
u(x,t=0)=\mbox{sin}[2\pi k (x-x_0)] e^{\f{(x-x_0)^2}{\sigma^2}}
\end{eqnarray}
where $x_0$ is the center of the wavepacket, initially located 
in the $x=0$. The parameter $\sigma$ is chosen such that tail 
of the wavepacket will be small enough inside the domain. In 
present study, we choose $\sigma=\sqrt{2}/10$. The computational 
domain is set to [-1,1] with periodic boundary conditions. If 
$c=1$, the configuration of the initial wavepacket repeats itself
every 2 time units. We consider a variety of wavenumbers $k$ to 
explore the resolution limit of the present CFOR-Hermite scheme. The 
mesh size is chosen as $1/N$ while $N=100$. Two time increments are 
selected to ensure that the time discretization error is negligible.
As our scheme is extremely accurate, it requires a very small time
increment to fully demonstrate its machine precision when $k=5$. 
Both the exact solution and numerical approximation at $t=100$ are 
plotted in FIG. \ref{wavepacket} for a comparison. Obviously, there is
no visual difference between them. Our results are listed in Tables 
\ref{advctab1} and \ref{advctab2} for a quantitative  evaluation.
For a large $k$ value, the wavepacket is highly oscillatory. It 
is a challenge for a low resolution scheme to resolve the advective 
wave. We explore the limitation of the present scheme by choosing
a variety of $k$ values. The performance of the CFOR-Hermite scheme 
can be well predicted by the Fourier analysis as shown in  FIG. 
\ref{advcfig1}.  It is noted that the Fourier transform of a 
Gaussian is still a Gaussian. The accuracy of the results is limited 
by the Fourier resolution of the Hermite differentiation, whose 
errors are indicated in FIG. \ref{advcfig1} for a number of 
different wavenumbers. For example, the case with $k=20$ is safely 
located to the left of the $\Delta\omega=10^{-10}$, which is the 
error of the differentiation. As a result, the CFOR error for the 
case of $k=20$ is of the order of $10^{-10}$ or less,  with  a small 
time increment. However, for the case of $k=30$,  its Gaussian tail 
is truncated by the limit $\Delta\omega=10^{-4}$. Consequently, the 
highest accuracy can not exceed this limitation, even with a small 
time increment.

The long time behavior of the CFOR-Hermite scheme for this problem 
is presented in Table \ref{advctab2}. We integrate the system
as long as 100 time units to investigate the reliability of the 
present approach for some large $k$ values ($k=20$ and 25).
The corresponding PPW values are 5 and 4 respectively. Both the 
long time stability and numerical resolution demonstrated in this 
problem are very difficult for most existing shock-capturing 
schemes to achieve.

{\sc Example 4.} {\it Advective 2D isentropic  vortex}

In this case, an isentropic vortex is introduced to a uniform 
mean flow, by small perturbations in velocity, density and the 
temperature. These perturbations involve all the three types 
of waves, i.e., vorticity, entropy and acoustic waves. The 
background flow and the perturbation parameters are given,
respectively by
\begin{eqnarray}
(u_\infty, v_\infty, p_\infty, T_\infty)  =  (1,1,1,1), \\
\end{eqnarray}
and 
\begin{eqnarray}
u'  =  -\f{\lambda}{2 \pi}(y-y_0) e^{\eta(1-r^2)}, \\ 
v'  =   \f{\lambda}{2 \pi}(x-x_0) e^{\eta(1-r^2)}, \\
T'  =  -\f{(\gamma-1) \lambda^2}{16 \eta \gamma \pi^2} 
       e^{2 \eta (1-r^2)},
\end{eqnarray}
where $r=\sqrt{(x-x_0)^2+(y-y_0)^2}$ is the distance to the 
vortex center, $\lambda=5$ is the strength of the vortex and 
$\eta$ is a parameter determining the gradient of the solution, 
and is unity in this study. By the isentropic relation, 
$p=\rho^\gamma$ and $T=p/\rho$, the perturbation in $\rho$ is 
required to be
\begin{eqnarray}
\rho=(T_{\infty} + T')^{1/(\gamma-1)} =
\l [ 1-\f{(\gamma-1) \lambda^2}{16 \eta \gamma \pi^2} 
    e^{2 \eta (1-r^2)}
\r ] ^{1/(\gamma-1)}. 
\label{density}
\end{eqnarray}
The governing equation is the compressible Euler equations with 
periodic boundary conditions\cite{eric01}. This problem was  
considered to measure the formal order and the stability in 
long-time integrations of many high order schemes\cite{eric01} 
because it is analytically solvable.

To be consistent with the literature\cite{eric01}, we chose a 
computational domain of $[0,10] \times [0,10]$ with the initial 
vortex center located at $(5,5)$. Our results  are tabulated in 
Tables \ref{vortextable1} and \ref{vortextable2}. Results of many 
standard schemes\cite{eric01} are also listed in these tables for 
a comparison. To be consistent with Ref. \cite{eric01} the error 
measures used in this problem are defined as 
\begin{eqnarray}
L_1 & = & {1\over (N+1)^2} 
  \sum_{i=0}^{N} \sum_{j=0}^{N} | f_{i,j}- \bar{f}_{i,j} |  \\
L_2 & = & {1\over (N+1)} \left[\sum_{i=0}^{N} \sum_{j=0}^{N}
               | f_{i,j}- \bar{f}_{i,j} |^2\right]^{1\over2},
\end{eqnarray}
where $f$ is the numerical result and $\bar{f}$ the exact solution.
Note that these definitions differ from the standard ones used 
in other problems. In Ref. \cite{eric01}, the CFL=0.5 was used 
and was nearly optimal. In our computation, the optimal time 
increment is much smaller because the nature of high resolution. 
We have set CFL=0.01 in some of our computations to release the 
full potential of the present scheme. It is noted that there is 
a dramatic gain in the accuracy as the mesh is refined from 
$40 \times 40$ to $80 \times 80$. The numerical order reaches
15 indicating that the proposed scheme is of extremely high order.
The present results obtained at the $80^2$ mesh are about $400$ 
times more accurate than those of other standard schemes obtained 
at $320^2$. Therefore, by using the present approach, a dramatical 
reduction of the mesh size and computing time can be gained for 
3D large scale computations.

The stability of the time integration is another useful test for
numerical schemes. We choose a mesh of $80\times80$ to test
the accuracy at long time computations with CFL=0.5. Although 
the time increment is not optimal for present scheme, practical 
computations usually require the time increment to be as large as
possible, especially for long time integrations. Numerical 
solutions are sampled at $t=2, 10, 50$ and 100. To prevent the
numerical errors from their nonlinear growth, the conjugate 
low-pass filter is activated to stabilize the computation. This 
problem differs from the advective sine-Gaussian wavepacket, where 
the equation is linear and the growth of the numerical error is 
fairly slow. The superior stability of the present CFOR scheme 
can be observed from Table \ref{vortextable3} and FIG. 
\ref{vortexfig1}. Even at $t=100$, the accuracy is still 
extremely high and the vortex core is well preserved.

{\sc Example 5.} {\it Shock/entropy wave interaction}

To further examine the resolution of the CFOR-Hermite scheme 
in the presence of shock waves, the interaction of shock/entropy 
waves is employed. In this case, small amplitude, low-frequency 
entropy waves are produced in the upstream of a shock. After 
the interaction with the shock, these waves are compressed in 
frequency and amplified in the amplitude. This is a standard test 
problem for high-order shock-capturing schemes and its governing 
equation is the compressible Euler system\cite{Shu1}. A linear 
analysis can be used to estimate the wave amplitude after the 
shock, which depends only on the shock strength and the amplitude 
before the shock. Therefore, the performance of a shock-capturing 
scheme can be evaluated by checking against the linear analysis. 
Low order methods perform poorly due to the oscillatory feature 
of the compressed waves. The ratio of frequencies before and after 
the shock depends only on the shock strength. The problem is defined 
in a domain of [0,5] and with following initial conditions\cite{Shu1}
\begin{eqnarray} \label{inicond}
 (\rho, u, p)  = \left \{
         \begin{array}{lcrl}
         (3.85714,                  & 2.629369, & 10.33333 ); &
                          x \le 0.5 \\
         (e^{- \epsilon \sin(\kappa x)}, & 0,        & 1.0 ); &
                          x > 0.5 \\
        \end{array} \right.,
\end{eqnarray}
where $\epsilon$ and $\kappa$ are the amplitude and wavenumber 
of the entropy wave before the shock. These initial conditions
represent a Mach 3 shock. We fix the amplitude of the wave 
$\epsilon$ before the shock to be 0.01. Accordingly the amplitude
after the shock is determined to be 0.08690716 by the linear 
analysis. In our previous study\cite{guwei}, the CFOR-RSK  
scheme was successfully applied to this problem. Here we focus 
on the resolution of the Hermite kernels. A variety of pre-shock 
wavenumbers $\kappa$ are considered for this purpose and the 
corresponding parameters are summarized in Table \ref{shocktab1}. 
Uniform meshes with grid spacing $1/N$ and $r=2.55$ are used for all 
the cases and their entropy waves are plotted in FIG. \ref{shockfig1}.

The entropy waves after the shock are well predicted in all the 
test cases, as indicated by the fully developed amplitude of 
entropy  waves after the shock. A low order scheme or a scheme 
with excessive  dissipation will lead to dramatical damping of the 
entropy waves. For $\kappa$=13, a mesh of $N=400$ is enough to 
well resolve the problem. When such a mesh is used in case of 
$\kappa=26$, the resolution is still very good. On a refined mesh 
$N=800$, the $\kappa=26$ case is better resolved as shown in FIG. 
\ref{shockfig1}(d). In our previous test using the CFOR-RSK scheme, 
the case of $\kappa=52$ was quite difficult to resolve on the mesh 
of  $N=800$. Such a difficulty is postponed in the present 
CFOR-Hermite scheme. The entropy waves of the case $\kappa=52$ 
with $N=800$ are plotted in FIG. \ref{shockfig1}(e), where the 
resolution is very good. Obviously, the good performance of 
the present CFOR-Hermite is attributed to its fast decay feature  
of the low-pass filter in its Fourier frequency response, see FIG.
\ref{compfig1}. When a finer mesh of  $N=1200$ is used for the 
case of $\kappa=52$, the quality of the solution near the shock 
front improves remarkably. Keeping the PPW being 5, a mesh of 
$N=1000$ is used for the case of $\kappa=65$ and results with 
high resolution are obtained, as shown in FIG. \ref{shockfig1}(g).
It can be noted that there is a slight over-dumping at the tail 
of the generated entropy waves and a remarkable improvement is 
achieved as the mesh is refined from $N=1000$ to $N=1200$, see 
FIG. \ref{shockfig1}(h). As the last example, the case of  
$\kappa=70$ is considered on the same mesh. Obviously, the 
CFOR-Hermite scheme yields a good solution for such a high 
wavenumber too. To our knowledge, such a high resolution has not 
been reported in the literature yet.

\section{Conclusion Remarks}

This paper introduces the conjugate filter oscillation reduction 
(CFOR) scheme with a Hermite kernel for the numerical simulation 
of both incompressible and compressible flows. The CFOR-Hermite 
scheme is based on a discrete singular convolution (DSC) algorithm,
which is a potential approach for the computer realization of 
singular convolutions. In particular, the DSC algorithm provides 
accurate approximations for spatial derivatives, which are high-pass 
filters form the point view of signal processing, in the smooth 
region of the solution. Therefore, the DSC algorithm makes it 
possible to resolve complex flow structures with a very small
ratio of point-per-wavelength (PPW) in fluid dynamical simulations. 
The essential idea of the CFOR scheme is that, a similar 
high-resolution low-pass filter which is conjugated to the high-pass 
filters will be activated to eliminate high-frequency errors  
whenever the spurious oscillation is generated in the flow for 
whatever reasons. The conjugate low-pass and high-pass filters 
are nearly optimal for shock-capturing and spurious oscillation 
suppression in the sense that they are generated from the same 
expression and consequently have similar order of regularity, 
effective frequency band and compact support.

An efficient and flexible DSC kernel plays a key role in the 
present CFOR scheme for resolving problems involving sharp 
gradients or shocks. Fourier analysis is conducted for both 
regularized Shannon kernel (RSK) and Hermite kernel. A comparison 
of both kernels reveals that the RSK has a wider effective 
frequency band for its high-pass filters,  whereas the Hermite 
kernel gives a better approximation to the idea low-pass filter. 
As a consequence, the DSC-RSK scheme achieves a better PPW 
ratio for flows without spurious oscillations, while the proposed
CFOR-Hermite scheme performs better for shock-capturing with 
natural high frequency oscillations. The effective frequency bands
of both kernels are much wider than those of most prevalent 
spatial discretization schemes.

Extensive numerical examples are employed to access the accuracy,
test the stability, demonstrate the usefulness and  most 
importantly, examine the resolution of the present CFOR-Hermite 
scheme. Our test cases cover four types of flow systems: 2D 
incompressible Euler systems, 1D linear advective equation, 2D 
compressible Euler systems and 1D shock/entropy wave interaction. 
Comparison is made between the RSK and Hermite kernel for the
2D incompressible flow. The DSC-RSK scheme performs better 
than the DSC-Hermite for smooth initial values. The capability 
of the CFOR-Hermite scheme for flow simulations is illustrated
by a variety of test examples. It is about $10^5$ times more 
accurate than some popular shock-capturing schemes in the 
advancement of a 2D isentropic vortex flow. Its advantage over 
the CFOR-RSK scheme is exemplified in predicting the shock/entropy 
wave interaction, where reliable results are attained at about 
5 PPW. Such a small ratio indicates that the present CFOR-Hermite
has a great potential for being used for large scale computations
with a very small mesh to achieve a given resolution.

\centerline{\bf Acknowledgment}

This work was supported in part by the National 
University of Singapore.


\centering
\begin{table}
\begin{tabular}{cccccccc}
\multirow{2}{1cm}{Error} & \multirow{2}{1cm}{Kernel} 
     & \multicolumn{6}{c}{$k$(PPW)} \\ \cline{3-8}
  &  & 1(32) &  2(16)  &  5(6.4) & 10(3.2) & 13(2.5) & 15(2.1)  \\ \hline
\multirow{2}{1cm}{$L_2$}      & Hermite
  & 6.63E-15 & 9.53E-15 & 2.45E-14 & 6.74E-13  & 1.01E-5 &  -  \\
  & RSK
  & 4.88E-15 & 3.55E-15 & 1.96E-14 & 1.47E-12 & 5.89E-11 & 1.55E-3 \\ \hline
\multirow{2}{1cm}{$L_\infty$} & Hermite
  & 2.78E-15 & 4.66E-15 & 1.86E-14 & 5.26E-13  & 4.79E-6 & -   \\
  & RSK
  & 2.33E-15 & 3.55E-15 & 1.64E-14 & 9.69E-13 & 4.19E-11 & 8.37E-4  \\
\end{tabular}
\vspace{2mm}
\caption{$L_2$ and $L_\infty$ errors of Taylor problem with 
         Hermite kernel and RSK at $t=2$. The PPW ratio in the
         bracket is calculated according to the wavenumber of 
         the pressure.}
         \label{taylortab1}

\newpage
\begin{tabular}{ccccccc}
\multirow{2}{1cm}{Time} & \multicolumn{6}{c}{$k$(PPW)} \\ \cline{2-7}
    & 5(20)    & 10(10)   & 15(6.7) & 20(5)   & 25(4)   & 30(3.3)  \\ \hline
2   & 2.00E-11 & 3.47E-10 & 2.26E-9 & 9.01E-9 & 3.34E-8 & 4.71E-5  \\ \hline
4   & 4.01E-11 & 6.95E-10 & 4.53E-9 & 1.80E-8 & 6.68E-8 & 9.41E-5  \\ \hline
6   & 6.01E-11 & 1.04E-9  & 6.79E-9 & 2.70E-8 & 1.00E-7 & 1.41E-4  \\ \hline
8   & 8.02E-11 & 1.39E-9  & 9.06E-9 & 3.60E-8 & 1.34E-7 & 1.88E-4  \\ \hline
10  & 1.00E-10 & 1.74E-9  & 1.13E-8 & 4.51E-8 & 1.67E-7 & 2.35E-4  \\
\end{tabular}
\vspace{2mm}
\caption{$L_1$ error for the advective sine-Gaussian 
         wavepacket, $\Delta t=10^{-4}$.} 
       \label{advctab1}
\vspace{1cm}
\begin{tabular}{ccccccc}
\multirow{2}{1cm}{Time} & \multicolumn{6}{c}{$k$(PPW)} \\ \cline{2-7}
    & 5(20)    & 10(10)    & 15(6.7)  & 20(5)    & 25(4)   & 30(3.3)  \\ \hline
2   & 1.17E-14 & 4.77E-14  & 4.23E-14 & 5.86E-12 & 2.21E-8 & 4.70E-5  \\ \hline
4   & 2.11E-14 & 8.86E-14  & 8.23E-14 & 1.17E-11 & 4.43E-8 & 9.41E-5  \\ \hline
6   & 2.46E-14 & 1.36E-13  & 1.11E-13 & 1.76E-11 & 6.64E-8 & 1.41E-4  \\ \hline
8   & 3.19E-14 & 1.79E-13  & 1.46E-13 & 2.35E-11 & 8.86E-8 & 1.88E-4  \\ \hline
10  & 4.01E-14 & 2.27E-13  & 1.73E-13 & 2.93E-11 & 1.11E-7 & 2.36E-4  \\
\end{tabular}
\vspace{2mm}
\caption{$L_1$ error for the advective sine-Gaussian 
         wavepacket, $\Delta t=5 \times 10^{-6}$.} 
 \label{advctab2}
\vspace{1cm}         
\begin{tabular}{ccccccc}
Time &   & 10 & 20 & 50 & 80 & 100 \\ \hline
\multirow{2}{1cm}{$L_1$} 
   & $k$=20 & 4.51E-8 & 9.01E-8 & 2.25E-7 & 3.60E-7 & 4.51E-7 \\ \cline{2-7}
   & $k$=25 & 1.67E-7 & 3.34E-7 & 8.35E-7 & 1.34E-6 & 1.67E-6 \\ \hline
\multirow{2}{1cm}{$L_\infty$} 
   & $k$=20 & 2.78E-7 & 5.56E-7 & 1.39E-6 & 2.22E-6 & 2.78E-6 \\ \cline{2-7}
   & $k$=25 & 1.51E-6 & 3.02E-6 & 7.55E-6 & 1.21E-5 & 1.51E-5 \\  
\end{tabular}
\vspace{2mm}
\caption{Errors for the advective sine-Gaussian wavepacket in long-time 
         integrations.} 
\label{advctab3}
\newpage      
\begin{tabular}{cccccccccccc}
N   &   & CFOR$^1$ & CFOR$^2$  & C4  & ENO  & MUSCL   
& WENO    & ENO$^{\it ACM}$  
& MUSCL$^{\it ACM}$ & WENO$^{\it ACM}$ \\ \hline
40    & error & 2.37E-5 & 6.45E-6 & 1.13E-3 & 1.28E-3 
      & 2.39E-3 & 9.39E-4 & 7.81E-4 & 1.29E-3 & 6.11E-4 \\ \hline
\multirow{2}{5mm}{80}  & error & 4.73E-9 & 2.79E-10 & 5.78E-5 & 2.08E-4 
      & 5.99E-4 & 7.07E-5 & 6.68E-5 & 2.79E-4 & 4.58E-4 \\ 
& order & 12.29 & 14.50 & 4.29 & 2.62 & 1.99 & 3.73 & 3.55 & 2.19 & 3.74  \\ \hline
\multirow{2}{5mm}{160} & error & 3.34E-10 & 3.76E-11 & 3.79E-6 & 3.01E-5 
        & 1.26E-4 & 2.46E-6 & 7.84E-6 & 5.31E-5 & 2.95E-6 \\ 
& order & 3.82 & 2.89 & 3.93 & 2.79 & 2.25 & 4.84 & 3.09 & 2.40 & 3.97  \\ \hline
\multirow{2}{5mm}{320} & error & 5.12E-11 & 3.20E-11 & 2.41E-7 & 4.07E-6 
        & 2.26E-5 & 8.52E-8 & 6.82E-7 & 8.61E-6 & 2.13E-7 \\ 
& order & 2.71 & 0.23 & 3.97 & 2.89 & 2.47 & 4.85 & 3.52 & 2.62 & 3.79  \\
\end{tabular} 
\vspace{0.5cm}    
\caption{
$L_1$ error for the density at $t=2$; 
The CFL number is 0.01 for CFOR$^2$ and 0.5
for all the other schemes;
C4: fourth-order accurate, conservative centered scheme;
ENO: third-order;
MUSCL: third-order;
WENO: fifth-order;
XXX$^{\it ACM}$ denotes that the C4 scheme is used as the basic scheme 
with the XXX being 
a characteristic-based filter, while Harten's artificial compression 
method is used as a sensor to indicate the local numerical dissipation,
see Ref. \protect \cite{eric01}. }
\label{vortextable1}
\vspace{1cm}
\begin{tabular}{cccccccccccc}
N   &  & CFOR$^1$ & CFOR$^2$ & C4  & ENO  & MUSCL   & WENO    & ENO$^{\it ACM}$  
& MUSCL$^{\it ACM}$ & WENO$^{\it ACM}$ \\ \hline
40    &  error & 4.35E-5 & 1.80E-5 & 2.92E-3 & 4.09E-3 
      & 8.29E-3 & 3.16E-3 & 2.47E-3 & 4.05E-3 & 2.08E-3 \\  \hline
\multirow{2}{5mm}{80}  & error & 1.41E-8 & 1.06E-9 & 1.90E-4 & 6.75E-4 
      & 2.26E-3 & 2.64E-4 & 2.08E-4 & 1.14E-3 & 1.48E-4 \\  
& order & 11.59 & 14.05 & 3.94 & 2.60 & 1.88 & 3.58 & 3.57 & 1.83 & 3.81 \\ \hline
\multirow{2}{5mm}{160} & error & 1.03E-9 & 4.73E-10 & 1.23E-5 & 8.69E-5 
      & 5.91E-4 & 1.10E-5 & 2.51E-5 & 3.12E-4 & 9.44E-6 \\ 
& order & 3.78 & 1.16 & 3.95 & 2.96 & 1.94 & 4.58 & 3.05 & 1.87 & 3.97 \\ \hline
\multirow{2}{5mm}{320} & error & 4.14E-10 & 4.08E-10 & 7.84E-7 & 1.33E-5 
      & 1.31E-4 & 2.93E-7 & 2.19E-6 & 6.07E-5 & 6.85E-7 \\ 
& order & 1.31 & 0.21 & 3.97 & 2.71 & 2.17 & 5.23 & 3.52 & 2.36 & 3.78 \\
\end{tabular} 
\vspace{0.5cm}
\caption{ $L_2$ error for the density at $t=2$. See Table 
\ref{vortextable1} for  captions.} 
\label{vortextable2}
\newpage
\begin{tabular}{ccccc}
Time  &  2  & 10  & 50  & 100              \\ \hline
$L_1$ &  4.73E-9 & 1.23E-8 & 4.58E-8 & 1.05E-7    \\ \hline
$L_2$ &  1.41E-8 & 3.64E-8 & 1.41E-7 & 3.17E-7    \\ 
\end{tabular} 
\vspace{0.5cm}
\caption{ Errors for the density at different times (CFL=0.5, N=80).}
\label{vortextable3}
\vspace{2cm}
\begin{tabular}{cccccccccc}
Case No.  &  1  & 2   & 3   & 4   & 5   & 6    & 7    &  8   & 9    \\ \hline
$k$       & 13  & 13  & 26  & 26  & 52  & 52   & 65   &  65  & 70   \\ \hline
N         & 400 & 800 & 400 & 800 & 800 & 1200 & 1000 & 1200 & 1200 \\ \hline
PPW       & 10  & 20  & 5   & 10  & 5   & 7.5  & 5    &    6 & 5.58 \\
\end{tabular} 
\vspace{5mm}
\caption{A summary of test cases for the shock/entropy wave interaction.} 
\label{shocktab1}        
\end{table}

\end{document}